\RequirePackage{fix-cm}
\documentclass{svjour3}                     
\smartqed  
\usepackage[utf8]{inputenc} 
\usepackage[T1]{fontenc}    
\usepackage{hyperref}       
\usepackage{url}            
\usepackage{booktabs}       
\usepackage{amsfonts}       
\usepackage{lipsum}
\usepackage{amsmath}
\usepackage{microtype}
\usepackage{graphicx}
\usepackage{subfigure}
\usepackage{booktabs} 
\usepackage{cite}
\usepackage{comment}
\usepackage{hyperref}

\newcommand{\norm}[1]{\left\lVert#1\right\rVert}

\DeclareMathOperator*{\argmin}{arg\,min}

\usepackage[framemethod=TikZ]{mdframed}
\mdfsetup{frametitlealignment=\center}
\begin{document}

\title{Factor-$\sqrt{2}$ Acceleration of Accelerated Gradient Methods}

\author{Chanwoo Park        \and
        Jisun Park          \and
        {Ernest K. Ryu}  
}


\institute{Chanwoo Park \at
              Department of Statistics, Seoul National University 
              \\
              \email{chanwoo.park@snu.ac.kr}           
           \and
           Jisun Park \at
              Department of Mathematical Sciences, Seoul National University\\
        \email{colleenp0515@snu.ac.kr}
           \and
              Ernest K. Ryu \at
              Department of Mathematical Sciences, Seoul National University\\
              \email{eryu@snu.ac.kr}
}

\date{Received: date / Accepted: date}

\maketitle

\begin{abstract}
The optimized gradient method (OGM) provides a factor-$\sqrt{2}$ speedup upon Nesterov's celebrated accelerated gradient method in the convex (but non-strongly convex) setup. However, this improved acceleration mechanism has not been well understood; prior analyses of OGM relied on a computer-assisted proof methodology, so the proofs were opaque for humans despite being verifiable and correct.  In this work, we present a new analysis of OGM based on a Lyapunov function and linear coupling. These analyses are developed and presented without the assistance of computers and are understandable by humans. Furthermore, we generalize OGM's acceleration mechanism and obtain a factor-$\sqrt{2}$ speedup in other setups: acceleration with a simpler rational stepsize, the strongly convex setup, and the mirror descent setup.
\end{abstract}

\section{Introduction}
\label{Introduction}

Nesterov's celebrated accelerated gradient method (AGM) solves the problem of finding the minimum of an $L$-smooth convex function with an ``optimal'' accelerated $\mathcal{O}(1/k^2)$ complexity  \cite{10029946121}.
Surprisingly, AGM turned out to be not exactly optimal, but optimal only up to a constant.
The optimized gradient method (OGM) has a factor-$2$ smaller (better) worst-case guarantee and thereby requires factor-$\sqrt{2}$ fewer iterations to guarantee the same accuracy \cite{drori2014performance, kim2016optimized}.

However, this remarkable discovery has not been well understood. OGM was originally obtained through a computer-assisted methodology based on the performance estimation problem (PEP). The resulting convergence analyses involve arduous but elementary calculations that are verifiable but arguably not understandable by humans.


\paragraph{Contribution.}
In this work, we present human-understandable analyses of OGM. First, we show that the improved acceleration mechanism of OGM can be understood and analyzed through an unconventional Lyapunov function. We then use this insight to propose a new method that obtains the factor-$\sqrt{2}$ speedup in the strongly convex setup. Finally, we present a human-understandable derivation of OGM based on refining the linear coupling analysis of Allen-Zhu and Orecchia \cite{allen2014linear}, and generalize OGM to the mirror descent setup.

As minor contributions, we analyze the primary and secondary sequences of OGM through a single unified analysis; 
to the best of our knowledge, prior works provide two separate convergence proofs for $x$- and $y$-sequences. 
Moreover, we present a unified class of accelerated methods containing AGM and OGM through the linear coupling analysis.

\subsection{Definitions and prior work}
For $L>0$, a differentiable convex function $f\colon\mathbb{R}^n\rightarrow\mathbb{R}$ is $L$-smooth with respect to a norm $\|\cdot\|$ if
\[
\|\nabla f(x)-\nabla f(y)\|_*\le L\|x-y\|\qquad
\forall\,x,y\in \mathbb{R}^n,
\]
where $\|\cdot\|_*$ denotes the dual norm.
A convex function $f\colon\mathbb{R}^n\rightarrow\mathbb{R}$ is $\mu$-strongly convex if $f(x)-(\mu/2)\|x\|^2$ is convex \cite{nesterov2003introductory, ryu2020LSCOMO}.

Throughout this paper, we consider the problem 
\begin{align*}
\begin{array}{ll}
\underset{x\in \mathbb{R}^n}{\mbox{minimize}}  &f(x)
  \end{array}
\end{align*}
and make the following assumptions on $f\colon\mathbb{R}^n\rightarrow\mathbb{R}$:
\begin{itemize}
\item [\hypertarget{A1}{(A1)}] $f$ is convex, differentiable, and $L$-smooth with respect to $\|\cdot\|$ and
\item [\hypertarget{A2}{(A2)}] $f$ has a minimizer (not necessarily unique).
\end{itemize}
We write $x_\star$ for a minimizer of $f$ and $f_\star=f(x_\star)$ for the optimal value.

\paragraph{Nesterov's AGM.}
Nesterov's AGM is
\begin{align*}
  y_{k+1} &= x_{k} - \frac{1}{L}\nabla{f(x_{k})}
\\
  x_{k+1} &= y_{k+1} + \frac{\theta_{k}-1}{\theta_{k+1}}(y_{k+1}- y_k),
\end{align*}
where $y_0=x_0$, $\theta_0 = 1$, and $\theta_{k+1}^2 - \theta_{k+1} =  \theta_k^2$ for $k=0,1,\dots$
\cite{10029946121}.
We can equivalently write AGM as
\begin{align*}
  y_{k+1} &= x_{k} - \frac{1}{L}\nabla{f(x_{k})}
\\
  z_{k+1} &= z_{k} - \frac{\theta_{k}}{L}\nabla{f(x_{k})} 
\\
  x_{k+1} &= \left(1- \frac{1}{\theta_{k+1}}\right)y_{k+1} + \frac{1}{\theta_{k+1}}z_{k+1}
\end{align*}
with $z_0=x_0$ \cite{nesterov2005smooth}.

AGM can be generalized to use the relaxed parameter requirement $\theta_{k+1}^2 - \theta_{k+1} \le  \theta_k^2$ on the positive sequence $\{ \theta_k \}_{k=0}^\infty$.
The choice $\theta_k=(k+2)/2$ is a common instance.

In the setup where $f$ is furthermore $\mu$-strongly convex, 
Nesterov's AGM for the strongly convex setup (SC-AGM) is
\begin{align*}
    y_{k+1} &= x_{k} - \frac{1}{L}\nabla{f(x_{k})}
    \\
    x_{k+1} &= y_{k+1} + \frac{\sqrt{\kappa}-1}{\sqrt{\kappa}+1}(y_{k+1}- y_k)   
\end{align*}
for $k=0,1,\dots$, where $\kappa=L/\mu$ and $ y_0=x_0$ \cite{nesterov2003introductory}.


\paragraph{Optimized gradient method.}
OGM is 
\begin{align*}
  y_{k+1} &= x_{k} - \frac{1}{L}\nabla{f(x_{k})}
\\
  x_{k+1} &= y_{k+1} + \frac{\theta_{k}-1}{\theta_{k+1}}(y_{k+1}- y_k) + \frac{\theta_k}{\theta_{k+1}}(y_{k+1} - x_{k})
\end{align*}
for $k=0,1,\dots$, where $ y_0=x_0$ and $\{\theta_k\}_{k=1}^\infty$ is the same as that of AGM \cite{drori2014performance, kim2016optimized}.
We refer to
$\frac{\theta_{k}-1}{\theta_{k+1}}(y_{k+1}- y_k)$
as the \emph{momentum term} and 
$\frac{\theta_k}{\theta_{k+1}}(y_{k+1} - x_{k})$
as the \emph{correction term}.
The added correction term is the difference between AGM and OGM.
We can equivalently write OGM as
\begin{align*}
  y_{k+1} &= x_{k} - \frac{1}{L}\nabla{f(x_{k})}
\\
  z_{k+1} &= z_{k} - \frac{2\theta_k}{L}\nabla{f(x_{k})} 
\\
  x_{k+1} &= \left(1-\frac{1}{\theta_{k+1}}\right)y_{k+1} + \frac{1}{\theta_{k+1}}z_{k+1},
\end{align*}
where $z_0=x_0$ \cite{kim2016optimized}.
The factor $2$ in $z_{k+1}$ is the difference compared to AGM. 

The $y_k$-sequence of OGM exhibits a rate faster than that of AGM by a factor of $\sqrt{2}$.
This rate was proved in 
\cite{kim2017convergence}, and we also state it in Corollary~\ref{cor:ogm-best-rate}.
To clarify, the guarantee on the function value is smaller (better) by a factor of $2$, and, combined with the  $\mathcal{O}(1/k^2)$ iteration dependence, this represents a factor-$\sqrt{2}$ reduction in the number of iterations necessary to reach a given accuracy.

Furthermore, OGM's original presentation \cite{drori2014performance, kim2016optimized} involves what we refer to as the \emph{last-step modification} on the secondary sequence
\begin{align*}
\tilde{x}_{k+1} &= y_{k+1} + \frac{\theta_{k}-1}{\varphi_{k+1}}(y_{k+1}- y_k) + \frac{\theta_k}{\varphi_{k+1}}(y_{k+1} - x_{k})\\
    &= \left(1-\frac{1}{\varphi_{k+1}}\right)y_{k+1} + \frac{1}{\varphi_{k+1}}z_{k+1} ,
\end{align*}
where $\varphi_{k}^2 - \varphi_{k} - 2\theta_{k-1}^2 = 0$. 
The $\tilde{x}_{k}$-sequence of OGM exhibits a rate slightly better than OGM's $y_k$-sequence and is in fact exactly optimal \cite{drori2017exact}.
This rate was proved in the original presentation of OGM \cite{drori2014performance, kim2016optimized}, and we also state it in Corollary~\ref{cor:ogm-best-rate-secondary}. 
In this work, we present the first variant of OGM for the strongly convex setup.

\paragraph{$\theta_k$-sequence asymptotic characterization.}
Throughout the exposition of this work, we will use the following asymptotic characterization:
if $\theta_0 = 1$ and $\theta_{k+1}^2 - \theta_{k+1} = \theta_{k}^2$ for $k = 0,1, \dots$, then 
\begin{equation}
    \theta_k = \frac{k+\zeta+1}{2} +  \frac{\log k}{4} +  o(1)
    \label{eq:theta-asymp}
\end{equation}
as $k\rightarrow\infty$, where $\zeta\approx 0.646$.
While we suspect this result may be known, we could not find it in any reference. 
Therefore, we formally state and prove \eqref{eq:theta-asymp} as Lemma~\ref{thm:theta-asymp} in the appendix.

\paragraph{Computer-assisted derivation and analysis of OGM.}
OGM was originally obtained through a computer-assisted methodology based on the performance estimation problem (PEP); it was first discovered numerically \cite{drori2014performance} and then its analytical form and convergence analysis was found \cite{kim2016optimized}.
The PEP methodology's key insight is to optimize over the class of fixed-step first-order gradient methods, with the objective being the convergence guarantee. Surprisingly, this problem is semidefinite programming- (SDP-) representable and has a tightness guarantee \cite{taylor2017smooth}.
OGM was re-discovered by using the PEP to find a greedy first-order method simplified with a ``subspace-search elimination procedure'' \cite{drori2020efficient}.

However, these prior analyses of OGM, generated by computers, are verifiable but arguably not understandable by humans.
Moreover, as the analyses rely on finding analytical solutions to the SDPs arising from the PEP, they are inaccessible to those unfamiliar with the methodology.

\paragraph{Lyapunov analysis of AGM.}
Nesterov's original 1983 paper established the celebrated $\mathcal{O}(1/k^2)$ rate using a Lyapunov analysis \cite{10029946121}.
Subsequent works \cite{nesterov2005smooth, auslender2006interior, nesterov2008accelerating, baes2009estimate, nesterov2012efficiency, nesterov2003introductory, li2020revisit, WibisonoE7351} analyzed AGM and its variants through the ``estimate sequence'' technique, which many consider to be less transparent than Lyapunov analyses.
In recent years, there has been a renewed interest in studying accelerated methods via Lyapunov analyses
\cite{v015a004, taylor2019stochastic, beck2009fast,ahn2020nesterov,  su2014differential, aujol2017optimal, aujol2019rates, aujol2019optimal}.
In this work, we present the first Lyapunov analysis of OGM.

\paragraph{Linear coupling analysis of AGM.}
The interpretation of AGM as a \emph{linear coupling} between gradient descent and mirror descent was presented in \cite{allen2014linear}.
Specifically, AGM can be written as
\begin{align*}
    y_{k+1} &= \argmin_y\left\{\langle \nabla f(x_k), y-x_k\rangle  +\frac{L}{2}\norm{y-x_k}^2 \right\}
    \\
    z_{k+1} &= \argmin_{y } \left\{ V_{z_k}(y) + \langle \alpha_{k+1} \nabla f(x_k), y-x_k \rangle \right\}
    \\
    x_{k+1} &= (1-\tau_{k+1}) y_{k+1} + \tau_{k+1} z_{k+1},
\end{align*}
where $V_z$ is a Bregman divergence. The $y_k$-update can be viewed as a gradient descent update and the $z_k$-update can be viewed as a mirror descent update.
Mirror descent \cite{nemirovsky1983problem} was originally presented as a method that maps the current point to a dual space, performs a gradient update, and maps the point back to the primal space.
An alternate proximal form of mirror descent (which we use) was presented in \cite{beck2003mirror}.
An alternate ``dual averaging'' interpretation of mirror descent as a method that constructs a lower bound of the function was presented in  \cite{nesterov2009primal}.
The key insight of linear coupling is to carefully interpolate between mirror descent and gradient descent to obtain AGM.

Linear coupling has been used to obtain and analyze many extensions of AGM \cite{allen2016even, allen2016variance, allen2016using, allen2017katyusha}, but whether the linear coupling argument itself can be further refined seems not to have been studied.
In this work, we show that refining the linear coupling analysis naturally leads to OGM.

\paragraph{Tight inequalities.}
We informally say an inequality is tight if it cannot be improved without further assumptions and formally if it satisfies the ``interpolation conditions'' \cite{taylor2017smooth}.
The recent literature on performance estimation problem focuses on using tight inequalities to obtain proofs that are provably cannot be improved \cite{ryu2020operator, de2020worst, taylor2017exact, gu2020tight, kim2019accelerated, lieder2020convergence, taylor2019stochastic}.

The tight inequality we use is
\[
f(y)\ge f(x)+\langle \nabla f(x),y-x\rangle+\frac{1}{2L}\|\nabla f(x)-\nabla f(y)\|^2_*
\]
for all $L$-smooth convex function $f$ and $x,y\in\mathbb{R}^n$.
The linear coupling analysis of AGM uses strictly weaker inequalities three times.
By refining the analysis by replacing the non-tight inequalities with tight ones, we obtain OGM.


\paragraph{Accelerated methods for smooth strongly convex minimization.}
For the problem setup of minimizing smooth strongly convex functions, Nesterov's SC-AGM \cite{nesterov2003introductory} achieves the convergence rate $\mathcal{O}\left(\exp\left(-k/\sqrt{\kappa}\right)\right)$.
Recently, the triple momentum method (TMM) \cite{van2017fastest} with an improved $\mathcal{O}\left(\exp\left(-2k/\sqrt{\kappa}\right)\right)$ rate was presented.
The SC-OGM method we present in this work has a rate of $\mathcal{O}\left(\exp\left(-\sqrt{2}k/\sqrt{\kappa}\right)\right)$, between the rates of SC-AGM and TMM, while still having a simple Lyapunov function-based proof.
For strongly convex convex \emph{quadratic} functions, heavy ball method exhibits the rate $\mathcal{O}\left(\exp\left(-4k\sqrt{\kappa}\right)\right)$ \cite{nesterov2003introductory} and OGM-q exhibits the rate $\mathcal{O}\left(\exp\left(-2\sqrt{2}k/\sqrt{\kappa}\right)\right)$ \cite{kim2018adaptive}.

\section{Lyapunov analysis of OGM}
\label{section:ogm}
In this section, we present a Lyapunov analysis of OGM. 
Our key insight is to use
\[
\left(f(x_{k}) - f_\star-\frac{1}{2L}\|\nabla f(x_k)\|^2\right),
\]
which is nonnegative due to $L$-smoothness,
instead of $\left(f(x_{k}) - f_\star\right)$ or $\left(f(y_{k}) - f_\star\right)$ in the construction of the Lyapunov function.
Throughout this section, $\|\cdot\|=\|\cdot\|_*$ denotes the Euclidean norm.

Based on this insight, we present:
(i) a more human-understandable analysis of OGM
(ii) a unified analysis of both the primary and secondary sequences of OGM that admits simpler $\theta_k$-choices.

\subsection{Nesterov's AGM}
Nesterov's AGM has the rate
\begin{align*}
    f(&y_{k}) - f_\star \leq \frac{L\norm{x_0 - x_\star}^2}{2\theta_{k-1}^2} =
    \frac{2L\norm{x_0 - x_\star}^2}{(k+\zeta)^2} - \frac{2L\norm{x_0 - x_\star}^2\log k}{(k+\zeta)^3} + o\left(\frac{1}{k^3}\right)
\end{align*}
for $k=0,1,\dots$.
This rate can be established through the following Lyapunov analysis \cite{10029946121}:
for $k=0,1,\dots$, define 
\begin{align*}
    U_k = \theta_{k-1}^2\left(f(y_{k}) - f_\star\right) + \frac{L}{2} \norm{z_{k} - x_\star}^2
\end{align*}
with $\theta_{-1}=0$ and show $U_k\le\cdots\le U_0$.
Conclude with
\[
\theta_{k-1}^2\left(f(y_{k}) - f_\star\right)\le U_k\le U_0=\frac{L}{2} \norm{x_0 - x_\star}^2.
\]

\subsection{Primary sequence analysis of OGM}
We now analyze OGM's convergence through an analogous Lyapunov analysis.
\begin{theorem}
\label{thm:OGM-main-thm}
Assume (\hyperlink{A1}{A1}) and (\hyperlink{A2}{A2}).
Let the positive sequence $\{\theta_k\}_{k=0}^\infty$ satisfy $\theta_0=1$ and $0 \leq \theta_{k+1}^2 - \theta_{k+1} \leq \theta_{k}^2$ for $k= 0,1,\dots$.
OGM's $y_k$-sequence exhibits the rate
\begin{align*}
f(y_{k}) -f_\star &\leq \frac{L\norm{x_0 - x_\star}^2}{4\theta_{k-1}^2}
\end{align*}
for $k=1,2,\dots$.

\end{theorem}
\begin{proof}
Set $\theta_{-1} = 0$ and $x_{-1} = x_0$. For $k=-1,0,1,\dots $, define
\begin{align*}
    U_k = &2\theta_{k}^2\left(f(x_k) - f_\star - \frac{1}{2L}\norm{\nabla f(x_k)}^2 \right) + \frac{L}{2}\norm{z_{k+1}- x_\star}^2.
\end{align*}
We can show that $\{U_k\}_{k=-1}^\infty$ is nonincreasing.
Using $f(y_{k}) \leq f(x_{k-1}) - \frac{1}{2L}\norm{\nabla f(x_{k-1})}^2$, which follows from $L$-smoothness, we conclude the rate with
\begin{align*}
2\theta_{k-1}^2\left(f(y_k) - f_\star \right) &\le
2\theta_{k-1}^2\left(f(x_{k-1}) - f_\star - \frac{1}{2L}\norm{\nabla f(x_{k-1})}^2 \right) \\
&\le U_{k-1}\le U_{-1}=\frac{L}{2}\norm{ z_{0} - x_\star}^2
\end{align*}
for $k =1,2,\dots$. Now we complete the proof by showing that $\{U_k\}_{k=-1}^\infty$ is nonincreasing.
For $k = -1,0,1, \dots$, we have
\begingroup
\allowdisplaybreaks
\begin{align*}
U&_{k} - U_{k+1} 
\\&= 2\theta_{k}^2\left(f(x_{k}) - f_\star - \frac{1}{2L}\norm{\nabla f(x_{k})}^2\right)-2\theta_{k+1}^2 \left(f(x_{k+1}) - f_\star - \frac{1}{2L}\norm{\nabla f(x_{k+1})}^2\right) 
\\
& \qquad+\frac{L}{2}\norm{z_{k+1} - x_\star}^2 - \frac{L}{2}\norm{ z_{k+2} - x_\star}^2
\\
&=2\theta_{k}^2\left(f(x_{k}) - f_\star - \frac{1}{2L}\norm{\nabla f(x_{k})}^2\right)-2\theta_{k+1}^2\left(f(x_{k+1}) - f_\star - \frac{1}{2L}\norm{\nabla f(x_{k+1})}^2\right)
\\
& \qquad-\langle 2\theta_{k+1} \nabla f(x_{k+1}), x_\star - z_{k+1} \rangle - \frac{2}{L}\theta_{k+1}^2 \norm{\nabla f(x_{k+1})}^2
\\
&=2\theta_{k}^2\left(f(x_{k}) - f_\star - \frac{1}{2L}\norm{\nabla f(x_{k})}^2\right)-2\theta_{k+1}^2\left(f(x_{k+1}) - f_\star + \frac{1}{2L}\norm{\nabla f(x_{k+1})}^2\right)
\\
& \qquad-\langle 2\theta_{k+1} \nabla f(x_{k+1}), x_\star - z_{k+1} \rangle
\\
& \geq 2(\theta_{k+1}^2- \theta_{k+1}) \left(f(x_{k}) - f_\star - \frac{1}{2L}\norm{\nabla f(x_{k})}^2\right)
\\
& \qquad-2\theta_{k+1}^2\left(f(x_{k+1}) - f_\star + \frac{1}{2L}\norm{\nabla f(x_{k+1})}^2\right)-\langle 2\theta_{k+1} \nabla f(x_{k+1}), x_\star - z_{k+1} \rangle
\\
& = 2(\theta_{k+1}^2- \theta_{k+1}) \left(f(x_{k}) - f_\star - \frac{1}{2L}\norm{\nabla f(x_{k})}^2 - f(x_{k+1}) + f_\star - \frac{1}{2L}\norm{\nabla f(x_{k+1})}^2\right) 
\\
& \qquad-2\theta_{k+1} \left(f(x_{k+1}) - f_\star + \frac{1}{2L}\norm{\nabla f(x_{k+1})}^2 \right) -\langle 2\theta_{k+1} \nabla f(x_{k+1}), x_\star - z_{k+1} \rangle
\\
& = 2(\theta_{k+1}^2- \theta_{k+1}) \left(f(x_{k}) - f(x_{k+1}) - \frac{1}{2L}\norm{\nabla f(x_{k})}^2 - \frac{1}{2L}\norm{\nabla f(x_{k+1})}^2\right) 
\\
& \qquad + 2\theta_{k+1} \left(f_\star - f(x_{k+1}) - \frac{1}{2L}\norm{\nabla f(x_{k+1})}^2 + \langle \nabla f(x_{k+1}), x_{k+1} - x_\star \rangle \right) 
\\
&\qquad \qquad + 2\theta_{k+1} \langle \nabla f(x_{k+1}), z_{k+1} - x_{k+1} \rangle
\\
& \geq 2(\theta_{k+1}^2- \theta_{k+1}) \left(f(x_{k}) - f(x_{k+1}) - \frac{1}{2L}\norm{\nabla f(x_{k})}^2 - \frac{1}{2L}\norm{\nabla f(x_{k+1})}^2\right) 
\\
&\qquad + 2\theta_{k+1} \langle \nabla f(x_{k+1}), z_{k+1} - x_{k+1} \rangle,
\end{align*}
\endgroup
where the inequalities follow from the cocoercivity of $f$.

Consider two separate cases $k = -1$ and $k = 0, 1, \dots$.
In case of $k = -1$, $\theta_{k+1}^2 - \theta_{k+1} = 1-1=0$ and $z_{k+1} - x_{k+1} = z_{0} - x_0 = 0$. The last formula becomes zero, so $U_{-1} - U_0 \geq 0$. 
For $k= 0,1,\dots$,
\begin{align*}
&2(\theta_{k+1}^2- \theta_{k+1}) \left(f(x_{k}) - f(x_{k+1}) - \frac{1}{2L}\norm{\nabla f(x_{k})}^2 - \frac{1}{2L}\norm{\nabla f(x_{k+1})}^2\right) 
\\
&\qquad + 2\theta_{k+1} \langle \nabla f(x_{k+1}), z_{k+1} - x_{k+1} \rangle
\\
& = 2(\theta_{k+1}^2- \theta_{k+1}) \left(f(x_{k}) - f(x_{k+1}) - \frac{1}{2L}\norm{\nabla f(x_{k})}^2 - \frac{1}{2L}\norm{\nabla f(x_{k+1})}^2\right) 
\\ & \qquad + 2\theta_{k+1}(\theta_{k+1} - 1) \langle \nabla f(x_{k+1}), x_{k+1} - x_{k}+ \frac{1}{L}\nabla f(x_{k}) \rangle
\\
&  = (2\theta_{k+1}^2 - 2\theta_{k+1})\biggl( f(x_{k})-f(x_{k+1}) -\frac{1}{2L}{\norm{\nabla f(x_{k}) - \nabla f(x_{k+1})}}^2 
\\
&\qquad \qquad \qquad \qquad \qquad + \langle \nabla f(x_{k+1}),x_{k+1}-x_{k}\rangle \biggr) \geq 0,
\end{align*}
where the inequalities follow from the cocoercivity of $f$.
\qed\end{proof}

As with AGM, the optimal $\{\theta_k\}_{k=0}^\infty$ is given by $\theta_{k+1}^2 - \theta_{k+1} = \theta_{k}^2$, which was used in the original presentation of OGM \cite{drori2014performance,kim2016optimized}.
\begin{corollary}
\label{cor:ogm-best-rate}
Under the setup of Theorem~\ref{thm:OGM-main-thm},
the choice $\theta_{k+1}^2 - \theta_{k+1} = \theta_{k}^2$ leads to the rate
\begin{align*}
f(&y_{k}) -f_\star \leq \frac{L\norm{x_0 - x_\star}^2}{4\theta_{k-1}^2}=
    \frac{L\norm{x_0 - x_\star}^2}{(k+\zeta)^2} - \frac{L\norm{x_0 - x_\star}^2 \log k}{(k+\zeta)^3} + o\left(\frac{1}{k^3}\right)
\end{align*}
for $k=1,2,\dots$.
\end{corollary}
\begin{proof}
This follows from Theorem \ref{thm:OGM-main-thm} and \eqref{eq:theta-asymp}.
\qed\end{proof}

The relaxed parameter requirement $0 \leq \theta_{k+1}^2 - \theta_{k+1} \le  \theta_{k}^2$ of Theorem~\ref{thm:OGM-main-thm} is reminiscent of the requirement for AGM.
We note that \cite{kim2018generalizing} had presented a generalized analysis with requirement $\theta_{k+1}^2 \leq \sum_{i=1}^{k+1} \theta_i$ based on the performance estimation problem methodology.

The relaxed parameter requirement allows us to use the simpler rational coefficients $\theta_k = (k+2)/2$. This leads to 
\begin{align*}
  y_{k+1} &= x_{k} - \frac{1}{L}\nabla{f(x_{k})}
\\
  x_{k+1} &= y_{k+1} + \frac{k}{k+3}(y_{k+1}- y_k) + \frac{k+2}{k+3}(y_{k+1} - x_{k}),
\end{align*}
which we call \emph{Simple-OGM}.

\begin{corollary}
Assume (\hyperlink{A1}{A1}) and (\hyperlink{A2}{A2}).
Simple-OGM's $y_k$-sequence exhibits the rate
\begin{align*}
f(y_{k}) -f_\star &\leq \frac{L\norm{x_0 - x_\star}^2}{(k+1)^2}
\end{align*}
for $k=1,2,\dots$.
\end{corollary}
\begin{proof}
This follows from Theorem~\ref{thm:OGM-main-thm}.
\qed\end{proof}

\subsection{Secondary sequence analysis of OGM}
We now analyze the convergence of OGM's secondary sequence with last-step modification through a unified Lyapunov analysis.

\begin{theorem}
\label{thm:OGM-main-thm-secondary}
Assume (\hyperlink{A1}{A1}) and (\hyperlink{A2}{A2}).
Let the positive sequence $\{\theta_k\}_{k=0}^\infty$ satisfy $\theta_0=1$, and $0 \leq \theta_{k+1}^2 - \theta_{k+1} \leq \theta_{k}^2$ for $k=0,1,\dots$.
Let the positive sequence $\{\varphi_k\}_{k=0}^\infty$ satisfy $0 \leq \varphi_{k}^2 - \varphi_{k} \leq 2\theta_{k-1}^2$ for $k=0,1,\dots$, where we define $\theta_{-1}=0$.
OGM's $\tilde{x}_k$-sequence, the secondary sequence with last-step modification, exhibits the rate
\begin{align*}
    f(\tilde{x}_{k}) - f_\star &\leq \frac{L\norm{x_0 - x_\star}^2}{2\varphi_k^2}
\end{align*}
 for $k=0,1,\dots$.
\end{theorem}
\begin{proof}
Let $\{U_{k}\}_{k=-1}^\infty$ be as defined in the proof of the Theorem~\ref{thm:OGM-main-thm}. Define $\{\tilde{U}_{k}\}_{k=0}^\infty$ as 
\begin{align*}
        \tilde{U}_{k} =& \varphi_{k}^2\left(f(\tilde{x}_{k}) - f_\star\right)+\frac{L}{2}\norm{ z_{k} - \frac{1}{L}\varphi_k \nabla f(\tilde{x}_{k}) - x_\star}^2.
\end{align*}
We can show that $\tilde{U}_k\le U_{k-1}$,
we conclude the rate with
\begin{align*}
\varphi_{k}^2\left(f(\tilde{x}_k) - f_\star \right)\le \tilde{U}_k\le U_{-1}=\frac{L}{2}\|x_0-x_\star\|^2
\end{align*}
for $k=0,1,\dots$. Now we complete the proof by showing that $\tilde{U}_k\le U_{k-1}$. For $k=0,1,\dots$, we have 
\begingroup
\allowdisplaybreaks
\begin{align*}
U&_{k-1} - \tilde{U}_{k} 
\\&= 2\theta_{k-1}^2\left(f(x_{k-1}) - f_\star - \frac{1}{2L}\norm{\nabla f(x_{k-1})}^2\right)-\varphi_{k}^2 \left(f(\tilde{x}_{k}) - f_\star\right) \\
& \qquad+\frac{L}{2}\norm{z_{k}-x_\star}^2 - \frac{L}{2}\norm{{z}_{k}- \frac{1}{L}\varphi_{k} \nabla f(\tilde{x}_{k})-x_\star}^2
\\&= 2\theta_{k-1}^2\left(f(x_{k-1}) - f_\star - \frac{1}{2L}\norm{\nabla f(x_{k-1})}^2\right)-\varphi_{k}^2 \left(f(\tilde{x}_{k}) - f_\star\right) \\
& \qquad-\langle \varphi_{k} \nabla f(\tilde{x}_{k}), x_\star - z_k \rangle - \frac{1}{2L}\varphi_{k}^2 \norm{\nabla f(\tilde{x}_{k})}^2
\\
&=2\theta_{k-1}^2\left(f(x_{k-1}) - f_\star - \frac{1}{2L}\norm{\nabla f(x_{k-1})}^2\right)
\\
& \qquad-\varphi_{k}^2\left(f(\tilde{x}_{k}) - f_\star + \frac{1}{2L}\norm{\nabla f(\tilde{x}_{k})}^2\right)-\langle \varphi_{k} \nabla f(\tilde{x}_{k}), x_\star - z_k \rangle
\\
& \geq (\varphi_{k}^2- \varphi_{k}) \left(f(x_{k-1}) - f_\star - \frac{1}{2L}\norm{\nabla f(x_{k-1})}^2\right)
\\
& \qquad-\varphi_{k}^2\left(f(\tilde{x}_{k}) - f_\star + \frac{1}{2L}\norm{\nabla f(\tilde{x}_{k})}^2\right)-\langle \varphi_{k} \nabla f(\tilde{x}_{k}), x_\star - z_k \rangle
\\
& = (\varphi_{k}^2- \varphi_{k}) \left(f(x_{k-1}) - f_\star - \frac{1}{2L}\norm{\nabla f(x_{k-1})}^2 - f(\tilde{x}_{k}) + f_\star - \frac{1}{2L}\norm{\nabla f(\tilde{x}_{k})}^2\right)
\\
& \qquad + \varphi_{k} \left(f_\star - f(\tilde{x}_{k}) - \frac{1}{2L}\norm{\nabla f(\tilde{x}_{k})}^2 + \langle \nabla f(\tilde{x}_{k}), \tilde{x}_{k} - x_\star \rangle \right)
\\
&\qquad \qquad + \langle \varphi_{k} \nabla f(\tilde{x}_{k}), z_k- \tilde{x}_{k} \rangle
\\
& \geq (\varphi_{k}^2- \varphi_{k}) \left(f(x_{k-1}) - f(\tilde{x}_{k}) - \frac{1}{2L}\norm{\nabla f(x_{k-1})}^2  - \frac{1}{2L}\norm{\nabla f(\tilde{x}_{k})}^2\right)
\\
& \qquad+ \langle \varphi_{k} \nabla f(\tilde{x}_{k}), z_k- \tilde{x}_{k} \rangle
\\
& = (\varphi_{k}^2- \varphi_{k}) \left(f(x_{k-1}) - f(\tilde{x}_{k}) - \frac{1}{2L}\norm{\nabla f(x_{k-1})}^2  - \frac{1}{2L}\norm{\nabla f(\tilde{x}_{k})}^2\right)
\\
& \qquad+ \varphi_{k}(\varphi_{k} - 1)  \langle \nabla f(\tilde{x}_{k}), \tilde{x}_{k}- x_{k-1} + \frac{1}{L} \nabla f(x_{k-1}) \rangle
\\
& = (\varphi_{k}^2- \varphi_{k}) \left( f(x_{k-1})-f(\tilde{x}_{k}) -\frac{1}{2L}{\norm{\nabla f(x_{k-1}) - \nabla f(\tilde{x}_{k})}}^2 + \langle \nabla f(\tilde{x}_{k}),\tilde{x}_{k}-x_{k-1}\rangle \right) 
\\
&\geq 0,
\end{align*}
\endgroup
where the inequalities follow from the cocoercivity of $f$.
\qed\end{proof}

\begin{corollary}
\label{cor:ogm-best-rate-secondary}
Under the setup of Theorem~\ref{thm:OGM-main-thm-secondary},
the choice $\theta_{k+1}^2 - \theta_{k+1} = \theta_{k}^2$ and $\varphi_{k}^2 - \varphi_{k}= 2\theta_{k-1}^2$ leads to the rate
\begin{align*}
&f(\tilde{x}_{k}) -f_\star \leq \frac{L\norm{x_0 - x_\star}^2}{2\varphi_k^2}=
    \frac{L\norm{x_0 - x_\star}^2}{(k+\zeta + 1/\sqrt{2})^2} - \frac{L\norm{x_0 - x_\star}^2\log k}{(k+\zeta + 1/\sqrt{2})^3} + o\left(\frac{1}{k^3}\right)
\end{align*}
for $k=0,1,\dots$.
\end{corollary}
\begin{proof}
This follows from \eqref{eq:theta-asymp}, which implies $\varphi_k = \frac{k+\zeta+\frac{1}{\sqrt{2}}}{\sqrt{2}}+ \frac{\sqrt{2}\log k}{4} +o(1)$, and Theorem~\ref{thm:OGM-main-thm-secondary}.
\qed\end{proof}

Simple-OGM with the last-step modification is
\begin{align*}
  y_{k+1} &= x_{k} - \frac{1}{L}\nabla{f(x_{k})}\\
  x_{k+1} &= y_{k+1} + \frac{k}{k+3}(y_{k+1}- y_k) + \frac{k+2}{k+3}(y_{k+1} - x_{k})\\
  \tilde{x}_{k+1}&= y_{k+1} + \frac{k}{\sqrt{2}(k+2) + 1}(y_{k+1} - y_{k}) + \frac{k+2}{\sqrt{2}(k+2) + 1}(y_{k+1} - x_k),
\end{align*}
where $ x_0=y_0$.

\begin{corollary}
Assume (\hyperlink{A1}{A1}) and (\hyperlink{A2}{A2}).
Simple-OGM's $\tilde{x}_k$-sequence, the secondary sequence with last-step modification, exhibits the rate
\begin{align*}
f(\tilde{x}_{k}) -f_\star &\leq \frac{L\norm{x_0 - x_\star}^2}{(k+1+ 1/\sqrt{2})^2}
\end{align*}
for $k=0, 1, \dots$.
\end{corollary}
\begin{proof}
Use Corollary \ref{cor:ogm-best-rate-secondary} with $\theta_k = \frac{k+2}{2}$ and $\varphi_k = \frac{k+1+\frac{1}{\sqrt{2}}}{\sqrt{2}} $.
\qed\end{proof}


\subsection{Discussion}
We clarify that the presented Lyapunov analysis is a novel contribution, while the results themselves are mostly known \cite{kim2016optimized,kim2017convergence,kim2018generalizing}.

We emphasize two key points.
First is the somewhat unusual construction of the Lyapunov function.
This key insight will be used in the following section to present a novel method for the strongly convex setup.

The second point we emphasize is that we present a unified analysis of the primary and last-step-modified secondary sequences using the Lyapunov functions $U_k$ and $\tilde{U_k}$.
Prior works on the two sequences of AGM and OGM rely on two separate analyses \cite{kim2016optimized, kim2017convergence}.

\section{Strongly convex OGM}
\label{section:sc-ogm}
In this section, we present strongly convex OGM (SC-OGM), a novel method that provides a factor-$\sqrt{2}$ improvement over Nesterov's SC-AGM.
The method and its analysis are obtained with following the key insight of Section~\ref{section:ogm}:
use the OGM-type correction term in the method and use
\[
\left(f(x_{k}) - f_\star-\frac{1}{2L}\|\nabla f(x_k)\|^2\right)
\]
in the construction of the Lyapunov function.
Throughout this section, $\|\cdot\|=\|\cdot\|_*$ denotes the Euclidean norm.

Based on this insight, we present:
(i) a novel method SC-OGM and
(ii) a unified analysis of both the primary and secondary sequences of SC-OGM.

\subsection{Nesterov's SC-AGM}
Further assume $f$ is $\mu$-strongly convex and write $\kappa = L/\mu $.
SC-AGM's convergence rate
\begin{align*}
f(y_{k}) - f_\star 
&\leq \left(1+\frac{1}{\sqrt{\kappa} -1}\right)^{-k} \frac{\mu + L}{2}\norm{x_0 - x_\star}^2 =\mathcal{O}\left(\exp\left(-\frac{k}{\sqrt{\kappa}}\right)\right)
\end{align*}
can be established through the following Lyapunov analysis \cite{v015a004}.
For $k=0,1,\dots$, define
\[
U_k = \left(1+\frac{1}{\sqrt{\kappa} -1}\right)^{k}\left(f(y_{k}) - f_\star + \frac{\mu}{2}\norm{z_{k} - x_\star}^2\right)
\]
with $z_{k} = (\sqrt{\kappa}+1)x_{k} - \sqrt{\kappa}y_k$ and show $U_k\le\cdots\le U_0 \le  \frac{\mu + L}{2}\norm{x_0 - x_\star}^2 $.

\subsection{Primary-sequence analysis of SC-OGM}
We newly propose SC-OGM:
\begin{align*}
    y_{k+1} &= x_{k} - \frac{1}{L}\nabla{f(x_{k})}
\\
    x_{k+1} &= y_{k+1} + \frac{1}{2\gamma + 1}(y_{k+1}-y_k) + \frac{1}{2\gamma + 1}(y_{k+1} - x_{k})
\end{align*}
for $k=0,1,\dots$, where $ y_0 = x_0$ and  $\gamma = \frac{\sqrt{8\kappa+1}+3}{2\kappa -2}$.

\begin{theorem}
\label{thm:SC-OGM-main-thm}
Assume (\hyperlink{A1}{A1}), (\hyperlink{A2}{A2}), and that $f$ is $\mu$-strongly convex.
SC-OGM's $y_k$-sequence exhibits the rate
\begin{align*}
f(y_{k}) - f_\star  &\leq (1+\gamma)^{-k+1}  \frac{\mu+2L}{2} \norm{x_0 - x_\star}^2=\mathcal{O}\left(\exp\left(-\frac{\sqrt{2}k}{\sqrt{\kappa}}\right)\right)
\end{align*}
for $k=1,2,\dots$.
\end{theorem}
\begin{proof}
For $k=0,1,\dots$, define
\[
z_{k} = \frac{2\gamma+1}{\gamma}x_{k} -\frac{\gamma +1}{\gamma}y_k
\]
and
\begin{align*}
    U_k =(1+\gamma)^k\biggl( f(x_k) - f_\star- &\frac{1}{2L}\norm{\nabla{f(x_k)}}^2 + \frac{\mu}{2}\norm{z_{k+1} - x_\star}^2\biggr).
\end{align*}
We can show that $\{U_k\}_{k=0}^\infty$ is nonincreasing and $U_0 \leq \frac{\mu+2L}{2}\norm{x_0 - x_\star}^2$. 
Using $f(y_{k}) \leq f(x_{k-1}) - \frac{1}{2L}\norm{\nabla f(x_{k-1})}^2$, which follows from $L$-smoothness, we conclude the rate with
\begin{align*}
(1+\gamma)^{k-1} \left( f(y_k) - f_\star \right) 
&\le(1+\gamma)^{k-1} \left( f(x_{k-1}) - f_\star - \frac{1}{2L} \norm{\nabla f(x_{k-1})}^2\right) \\
&\le U_{k-1}\le U_{0} \leq \frac{\mu+2L}{2}\norm{x_0 - x_\star}^2
\end{align*}
for $k=1,2,\dots$. Now we complete the proof by showing $U_0 \leq \frac{\mu+2L}{2}\norm{x_0 - x_\star}^2$, showing some relationships between $x_k$ and $z_k$, and showing that $\{U_k\}_{k=0}^\infty$ is nonincreasing. 

Firstly, we have
\begingroup
\allowdisplaybreaks
\begin{align*}
    U_0 &= f(x_0) -f_\star - \frac{1}{2L} \norm{\nabla f(x_0)}^2 + \frac{\mu}{2}\norm{z_1 - x_\star}^2
    \\
    &= f(x_0) -f_\star - \frac{1}{2L} \norm{\nabla f(x_0)}^2 + \frac{\mu}{2}\norm{x_0 - \frac{1}{L}\frac{\gamma +2}{\gamma}\nabla f(x_0) - x_\star}^2
    \\
    &= f(x_0) - f_\star + \frac{1}{2L}\frac{1}{\gamma + 1} \norm{\nabla f(x_0)}^2 - \frac{\gamma}{1+\gamma}\langle \nabla f(x_0), x_0 - x_\star \rangle + \frac{\mu}{2}\norm{x_0 - x_\star}^2 
    \\
    &\leq f(x_0) - f_\star + \frac{1}{2L} \norm{\nabla f(x_0)}^2 - \frac{\gamma}{1+\gamma}\langle \nabla f(x_0), x_0 - x_\star \rangle + \frac{\mu}{2}\norm{x_0 - x_\star}^2 
    \\
    &\leq \frac{1}{\gamma + 1}(f(x_0) - f_\star) + \frac{1}{2L}\frac{1}{1+\gamma} \norm{\nabla f(x_0)}^2+ \frac{\mu}{2}\norm{x_0 - x_\star}^2  
    \\
    &\leq \frac{2}{1+\gamma}(f(x_0) - f_\star) + \frac{\mu}{2}\norm{x_0 - x_\star}^2  
    \\
    &\leq \left(L + \frac{\mu}{2}\right)\|x_0 - x_\star\|^2.
\end{align*}
\endgroup

Second, Let $X_k = x_k - x_\star$ and $Z_k = z_k - x_\star$, for $k = 0, 1, \dots$. We will prove 
\begin{align}
(x_{k+1} - x_{k}) + \frac{1}{L}\nabla{f(x_{k})} + \gamma X_{k+1} &= \frac{1}{1 + \gamma}(\gamma Z_{k+1} + \gamma ^2 X_{k+1}) \label{eq:scogm-xkzk1}
\\
Z_{k+1} &= \frac{1}{\gamma +1} Z_k +\frac{\gamma}{\gamma +1} X_{k} - \frac{1}{L}\frac{\gamma + 2}{\gamma} \nabla{f(x_{k})} \label{eq:scogm-xkzk2}
\end{align}
for $k = 0,1, \dots$.

Plug $y_{k+1} = x_{k} - \frac{1}{L}\nabla f(x_{k})$ in the definition of $z_{k+1}$. Then we obtain (\ref{eq:scogm-xkzk1}). 
\\
For (\ref{eq:scogm-xkzk2}), from definition of $z_k$ and $z_{k+1}$
\begin{align*}
z_{k+1} &=   \frac{2\gamma +1}{\gamma}x_{k+1} - \frac{\gamma + 1}{\gamma}x_{k} +\frac{1}{L}\frac{1+\gamma}{\gamma} \nabla{f(x_{k})}
\\
z_{k} &=   \frac{2\gamma +1}{\gamma}x_{k} - \frac{\gamma + 1}{\gamma}x_{k-1} +\frac{1}{L}\frac{1+\gamma}{\gamma} \nabla{f(x_{k-1})}
\end{align*}
and definition of $x_k$, we have
\begin{align*}
    x_{k+1} &= \frac{2\gamma + 2}{2\gamma + 1 }y_{k+1}-\frac{1}{2\gamma + 1 } y_{k}-\frac{1}{L}\frac{1}{2\gamma+1}\nabla{f(x_{k})}
    \\
    &=\frac{2\gamma + 2}{2\gamma + 1 }x_{k}-\frac{1}{2\gamma + 1 }x_{k-1}-\frac{1}{L}\frac{2\gamma + 3}{2\gamma +1} \nabla{f(x_{k})}+\frac{1}{L}\frac{1}{2\gamma +1}\nabla{f(x_{k-1})}.
\end{align*}
Therefore,
\begingroup
\allowdisplaybreaks
\begin{align*}
    z_{k+1} - \frac{1}{\gamma+1}z_{k} &= \frac{2\gamma +1}{\gamma}x_{k+1} - \frac{\gamma + 1}{\gamma}x_{k} +\frac{1}{L}\frac{1+\gamma}{\gamma} \nabla{f(x_{k})} 
    \\
    &\qquad - \frac{1}{\gamma+1}\left(\frac{2\gamma +1}{\gamma}x_{k} - \frac{\gamma + 1}{\gamma}x_{k-1} +\frac{1}{L}\frac{1+\gamma}{\gamma} \nabla{f(x_{k-1})} \right)
    \\
    &=\frac{2\gamma +1}{\gamma}x_{k+1} - \frac{\gamma^2 + 4\gamma + 2 }{\gamma(\gamma+1)}x_{k} + \frac{1}{\gamma}x_{k-1} +\frac{1}{L}\frac{1+\gamma}{\gamma}\nabla{f(x_{k})}
    \\
    &\qquad -\frac{1}{L}\frac{1}{\gamma} \nabla{f(x_{k-1})} 
    \\
    &=\frac{2\gamma +1}{\gamma}\biggl(\frac{2\gamma + 2}{2\gamma + 1 }x_{k}-\frac{1}{2\gamma + 1 }x_{k-1}-\frac{1}{L}\frac{2\gamma + 3}{2\gamma +1} \nabla{f(x_{k})}
    \\
    &\qquad +\frac{1}{L}\frac{1}{2\gamma +1}\nabla{f(x_{k-1})}\biggr) - \frac{\gamma^2 + 4\gamma + 2 }{\gamma(\gamma+1)}x_{k} + \frac{1}{\gamma}x_{k-1}
    \\
    &\qquad \qquad +\frac{1}{L}\frac{1+\gamma}{\gamma}\nabla{f(x_{k})} -\frac{1}{L}\frac{1}{\gamma} \nabla{f(x_{k-1})} 
    \\
    &= \frac{\gamma}{\gamma + 1}x_{k} - \frac{1}{L}\frac{\gamma + 2}{\gamma} \nabla{f(x_{k})}
\end{align*}
\endgroup
so we obtained (\ref{eq:scogm-xkzk2}).

Lastly, we will show that $\{U_k\}_{k=0}^\infty$ is nonincreasing. It suffices to show that for $k=0,1, \dots $,
$$ (1+\gamma)^{-k} (U_k - U_{k+1}) \geq 0 $$ 
which is equivalent to showing
\begin{align*}
&\left((f(x_{k}) - f_\star- \frac{1}{2L}\norm{\nabla{f(x_{k})}}^2) - (1+\gamma)(f(x_{k+1}) - f_\star-\frac{1}{2L}\norm{\nabla{f(x_{k+1})}}^2)\right)
\\
&\qquad \qquad+ \frac{\mu}{2}\left( \norm{z_{k+1} - x_\star}^2 - (1+\gamma)\norm{z_{k+2} - x_\star}^2 \right) \geq 0.
\end{align*}

By $L$-smoothness of $f$, we have 
$$ f(x_{k+1}) - f(x_{k}) \leq -\frac{1}{2L}\norm{\nabla{f(x_{k+1})-\nabla{f(x_{k})}}}^2 +\langle \nabla{f(x_{k+1})}, x_{k+1} - x_{k} \rangle $$
and from strong convexity, $$f(x_{k+1}) - f_\star \leq \langle\nabla{f(x_{k+1})}, x_{k+1} - x_\star\rangle  - \frac{\mu}{2} \norm{x_{k+1} - x_\star}^2.$$
For $k = 0,1, \dots$, using above two inequalities, \eqref{eq:scogm-xkzk1}, and \eqref{eq:scogm-xkzk2},

\begin{align*}
    \biggl(&f(x_{k}) - f_\star- \frac{1}{2L}\norm{\nabla{f(x_{k})}}^2\biggr) - (1+\gamma)(f(x_{k+1}) - f_\star-\frac{1}{2L}\norm{\nabla{f(x_{k+1})}}^2)
    \\
    &=(f(x_{k}) - f(x_{k+1})) - \gamma(f(x_{k+1}) - f_\star) + \frac{1+\gamma}{2L}\norm{\nabla{f(x_{k+1})}}^2 - \frac{1}{2L}\norm{\nabla{f(x_{k})}}^2
    \\
    &\geq \left(\frac{1}{2L}\norm{\nabla{f(x_{k+1})-\nabla{f(x_{k})}}}^2 +\langle \nabla{f(x_{k+1})}, x_{k} - x_{k+1} \rangle\right) 
    \\
    &\qquad -  \gamma \left(\langle\nabla{f(x_{k+1})}, x_{k+1} - x_\star\rangle - \frac{\mu}{2} \norm{x_{k+1} - x_\star}^2\right) 
    \\
    &\qquad \qquad + \frac{1+\gamma}{2L}\norm{\nabla{f(x_{k+1})}}^2 -\frac{1}{2L}\norm{\nabla{f(x_{k})}}^2
    \\
    &= \langle \nabla{f(x_{k+1})}, -\frac{1}{L}\nabla{f(x_{k})} - x_{k+1} + x_{k} - \gamma(x_{k+1} - x_\star) \rangle
    \\
    &\qquad + \frac{2+\gamma}{2L}\norm{\nabla{f(x_{k+1})}}^2 + \frac{\mu \gamma}{2}\norm{x_{k+1} - x_\star}^2 
    \\
    &= \langle \nabla{f(x_{k+1})}, -\frac{1}{1+\gamma}(\gamma Z_{k+1} + \gamma^2 X_{k+1})  \rangle
    \\
    &\qquad + \frac{2+\gamma}{2L}\norm{\nabla{f(x_{k+1})}}^2 + \frac{\mu \gamma}{2}\norm{x_{k+1} - x_\star}^2. 
\end{align*}

In addition,
\begin{align*}
\frac{\mu}{2}&\left((1+\gamma)\norm{Z_{k+2}}^2 - \norm{Z_{k+1}}^2\right) 
\\
&=\frac{\mu}{2}\left((1+\gamma)\norm{\frac{1}{1+\gamma}Z_{k+1} + \frac{\gamma}{1+\gamma}X_{k+1} - \frac{1}{L} \frac{2+\gamma}{\gamma}\nabla{f(x_{k+1})}}^2 - \norm{Z_{k+1}}^2 \right) 
\\
&= \frac{\mu}{2}(-\frac{\gamma}{1+\gamma}\norm{ Z_{k+1}}^2 + \frac{\gamma^2}{1+\gamma} \norm{X_{k+1}}^2 + (1+\gamma)\frac{1}{L^2}\frac{(2+\gamma)^2}{\gamma^2}\norm{\nabla{f(x_{k+1})}}^2 
\\&\qquad +2 \frac{\gamma}{1+\gamma}\langle Z_{k+1}, X_{k+1} \rangle  - 2\frac{2+\gamma}{L\gamma}\langle \nabla f(x_{k+1}), Z_{k+1} \rangle 
\\
&\qquad \qquad - 2\frac{2 + \gamma}{L}\langle \nabla f(x_{k+1}), X_{k+1} \rangle ).
\end{align*}
Since 
$$\mu \frac{2+\gamma}{L\gamma^2} = \frac{1}{1+\gamma}, $$
we can telescope concerned $\nabla f(x_{k+1})$'s inner product in $U_k - U_{k+1}$.

For $k =0,1, \dots$, we have
\begin{align*}
    (1+\gamma)^{-k} & (U_k - U_{k+1})
    \\
    &\geq   
    \frac{2+\gamma}{2L}\norm{\nabla{f(x_{k+1})}}^2 +  \frac{\mu \gamma}{2}\norm{X_{k+1}}^2   
    \\
    & \qquad - \frac{\mu}{2}\biggl(-\frac{\gamma}{1+\gamma}\norm{ Z_{k+1}}^2 + \frac{\gamma^2}{1+\gamma} \norm{X_{k+1}}^2 
    \\
    &\qquad \qquad + (1+\gamma)\frac{1}{L^2}\frac{(2+\gamma)^2}{\gamma^2}\norm{\nabla{f(x_{k+1})}}^2  +2 \frac{\gamma}{1+\gamma}\langle Z_{k+1}, X_{k+1} \rangle  \biggr)
    \\
    &= -\frac{\mu}{2} \left( - \frac{\gamma}{1+\gamma}\norm{X_{k+1}}^2 - \frac{\gamma}{1+\gamma}\norm{Z_{k+1}}^2  +2 \frac{\gamma}{1+\gamma}\langle Z_{k+1}, X_{k+1} \rangle\right)
    \\
    &=\frac{\mu}{2}\frac{\gamma}{1+\gamma}\norm{Z_{k+1} - X_{k+1}}^2 \geq 0.
\end{align*}
\qed\end{proof}

\subsection{Secondary sequence analysis}
We now analyze the convergence of SC-OGM's secondary sequence with a unified Lyapunov analysis.
We note that SC-OGM does not require the last-step modification, unlike the non-strongly convex counterpart.

\begin{theorem}
\label{thm:SC-OGM-main-thm-secondary}
Assume (\hyperlink{A1}{A1}), (\hyperlink{A2}{A2}), and that $f$ is $\mu$-strongly convex.
SC-OGM's $x_k$-sequence, the secondary sequence without last-step modification, exhibits the rate
\begin{align*}
f(x_{k}) - f_\star  \leq \frac{(1+\gamma)^{-k+2}}{2\gamma}  \left(\frac{\mu+2L}{2} \norm{x_0 - x_\star}^2\right)
\end{align*}
for $k=1,2, \dots$.
\end{theorem}
\begin{proof}
Let $\{z_k\}_{k=0}^\infty$ and $\{U_{k}\}_{k=0}^\infty$ be defined as in the proof of the Theorem~\ref{thm:SC-OGM-main-thm}.
For $k=0,1,\dots$, define
\begin{align*}
    \tilde{U}_k = (1+&\gamma)^{k-1}\biggl( \frac{2\gamma}{1+\gamma}(f(x_k) - f_\star) +  \frac{\mu}{2}\norm{z_k - \left(\frac{\gamma + 2}{\gamma}\right)\frac{1}{L}\nabla f(x_k) - x_\star}^2\biggr)
\end{align*}
We can show that $\tilde{U}_k\le U_{k-1}$. We conclude the rate with
\begin{align*}
(1+\gamma)^{k-1}\frac{2\gamma}{1+\gamma}\left( f(x_k) - f_\star \right) \le \tilde{U}_{k}\le U_{0} \leq \frac{\mu+2L}{2}\norm{x_0 - x_\star}^2
\end{align*}
for $k=1,2,\dots$. Now we complete the proof by showing that $\tilde{U}_k\le U_{k-1}$. Note that $\frac{\gamma + 1}{\gamma}\left((x_k - x_{k-1}) + \frac{1}{L}\nabla{f(x_{k-1})}\right) = (Z_k - X_k) $. Then we have
\begingroup
\allowdisplaybreaks
\begin{align*}
&\left(f(x_{k-1}) - f_\star - \frac{1}{2L}\norm{\nabla f(x_{k-1})}^2\right)-\frac{2\gamma}{1+\gamma}\left(f(x_k) - f_\star\right) 
\\ 
&\quad + \frac{L\gamma^2}{2(1+\gamma)(2+\gamma)}\norm{z_k - x_\star}^2 - \frac{L\gamma^2}{2(1+\gamma)(2+\gamma)}\norm{z_k - \left(\frac{\gamma + 2}{\gamma}\right)\frac{1}{L}\nabla f(x_k) - x_\star}^2
\\
&=\left(f(x_{k-1}) - f_\star - \frac{1}{2L}\norm{\nabla f(x_{k-1})}^2\right)-\frac{2\gamma}{1+\gamma}\left(f(x_k) - f_\star\right)  
\\
&\qquad + \frac{\gamma}{1+\gamma}\langle Z_k, \nabla{f(x_k)} \rangle - \frac{1}{2L} \frac{2 + \gamma}{1+\gamma}\norm{\nabla{f(x_k)}}^2 
\\
&=\left(f(x_{k-1}) - f_\star - \frac{1}{2L}\norm{\nabla f(x_{k-1})}^2\right)-\frac{2\gamma}{1+\gamma}\left(f(x_k) - f_\star\right)
\\
& \qquad  + \frac{\gamma}{1+\gamma} \left\langle \frac{\gamma + 1}{\gamma}\left((x_k - x_{k-1}) + \frac{1}{L}\nabla{f(x_{k-1})}\right) + X_k, \nabla{f(x_k)} \right\rangle 
\\
&\qquad \qquad - \frac{1}{2L} \frac{2 + \gamma}{1+\gamma}\norm{\nabla{f(x_k)}}^2 
\\
&=\left(f(x_{k-1}) - f_\star - \frac{1}{2L}\norm{\nabla f(x_{k-1})}^2\right)-\frac{2\gamma}{1+\gamma}\left(f(x_k) - f_\star\right) 
\\
& \qquad  + \langle x_k - x_{k-1}, \nabla{f(x_k)} \rangle + \frac{1}{L} \langle \nabla{f(x_{k-1})}, \nabla f(x_k) \rangle  + \frac{\gamma}{1+\gamma}\langle X_k, \nabla{f(x_k)} \rangle 
\\
&\qquad \qquad - \frac{1}{2L} \frac{2 + \gamma}{1+\gamma}\norm{\nabla{f(x_k)}}^2 
\\
&=\left( f(x_{k-1})-f(x_k) -\frac{1}{2L}{\norm{\nabla f(x_{k-1}) - \nabla f(x_k)}}^2 + \langle \nabla f(x_k),x_k-x_{k-1}\rangle \right)  
\\
&\qquad + \frac{1}{2L}\frac{\gamma}{1+\gamma}\norm{\nabla f(x_k)}^2+ \frac{1}{1 + \gamma} \left(f(x_{k}) - f_\star - \frac{1}{2L}\norm{\nabla f(x_{k})}^2\right) 
\\
& \qquad  \qquad  + \frac{\gamma}{1 + \gamma} \left(f_\star - f(x_k) - \frac{1}{2L} \norm{\nabla f(x_k)}^2 + \langle X_k, \nabla f(x_k) \rangle\right)
\\
&\geq 0.
\end{align*}
\endgroup
Since $\frac{L\gamma^2}{2(1+\gamma)(2+\gamma)} = \frac{\mu}{2}$, above inequality indicates that 
\begin{align*}
\biggl(f(x_{k-1}) - f_\star - &\frac{1}{2L}\norm{\nabla f(x_{k-1})}^2\biggr) + \frac{\mu}{2}\norm{z_k - x_\star}^2
\\&\geq \frac{2\gamma}{1+\gamma}(f(x_k) - f_\star) +  \frac{\mu}{2}\norm{z_k - \left(\frac{\gamma + 2}{\gamma}\right)\frac{1}{L}\nabla f(x_k) - x_\star}^2. 
\end{align*}
\qed\end{proof}

\subsection{Discussion}
The factor-$\sqrt{2}$ improvement of SC-OGM over SC-AGM is consistent with  the factor-$\sqrt{2}$ improvement of OGM over AGM.
AGM and OGM share the same momentum term while OGM has the additional ``correction term''.
In contrast, the momentum coefficients differ in the strongly convex case: SC-AGM has
\[
\frac{\sqrt{\kappa}-1}{\sqrt{\kappa}+1}=1-\frac{2}{\sqrt{\kappa}} + 
\mathcal{O}\left(\frac{1}{\kappa}\right)
\]
while SC-OGM has
\[
    \frac{1}{2\gamma + 1}=1-\frac{2\sqrt{2}}{\sqrt{\kappa}}+ \mathcal{O}\left(\frac{1}{\kappa}\right).
\]
Of course, SC-OGM also has the correction term, which is essential in the analysis.

Again, we point out that we analyze both the primary and secondary sequences in a single streamlined proof using the Lyapunov functions $U_k$ and $\tilde{U_k}$.

We also point out that \cite{kim2018adaptive} presented a variant of OGM named OGM-q for smooth strongly convex \emph{quadratic} functions.
In contrast, SC-OGM applies to the broader class of smooth strongly convex functions.


\section{Linear coupling analysis}
\label{section:4}
While the Lyapunov analyses of Sections~\ref{section:ogm} and \ref{section:sc-ogm} do provide insight into the acceleration mechanism of OGM,
they do not shed light onto the \emph{provenance} of the method.
Originally, OGM was generated through a computer-assisted proof methodology as the exactly optimal first-order method, but this approach is arguably opaque to humans.

In this section, we present a human-understandable \emph{derivation} of OGM based on linear coupling.
Specifically, we obtain OGM by refining the linear coupling analysis of Allen-Zhu and Orecchia \cite{allen2014linear} through replacing the use of non-tight inequalities with tight inequalities.

We specifically provide:
(i) a natural (and non-computer assisted) derivation of OGM, 
(ii) a generalization of OGM to the mirror descent setup, and
(iii) a unification of AGM and OGM.
We moreover provide (iv) a generalization of SC-OGM to the mirror descent setup in the appendix, in Section \ref{sec:lc-sc-ogm}.

\paragraph{Assumption and notation.}
In this section, assume
\begin{itemize}
    \item [\hypertarget{A3}{(A3)}] $\|\cdot\|=\sqrt{x^TQx}$ is a quadratic norm, where $Q$ is a symmetric positive definite matrix.
\end{itemize}
Assumption (\hyperlink{A1}{A1}) is to be interpreted as $L$-smoothness with respect to norm $\|\cdot\|$.
Write $\|\cdot\|_*=x^TQ^{-1}x$ for the dual norm of $\|\cdot\|$.
However, $\langle\cdot,\cdot\rangle$ is the standard Euclidean inner product (unrelated to $Q$).
Let $w\colon\mathbb{R}^n\rightarrow\mathbb{R}$ be a ``distance generating function'' that is differentiable and $1$-strongly convex with respect to $\|\cdot\|$,
and let
\[
V_x(y)=w(y)-\langle \nabla w(x),y-x\rangle-w(x)\qquad
\forall x,y\in\mathbb{R}^n
\]
be the Bregman divergence generated by $w$.

\subsection{Linear coupling analysis of AGM}
We briefly outline the linear coupling analysis of AGM presented in \cite{allen2014linear} and point out where the analysis can be refined.

Consider the problem of minimizing $f$ under assumptions (\hyperlink{A1}{A1}), (\hyperlink{A2}{A2}), and (\hyperlink{A3}{A3}). 
The linear coupling method is
\begin{align*}
    y_{k+1} &= x_k-L^{-1}Q^{-1}\nabla f(x_k)\tag{LC}
    \label{eq:linear-coupling}
    \\
    z_{k+1} &= \argmin_{y\in\mathbb{R}^n}\left\{V_{z_k}(y)+\langle \alpha_{k+1}\nabla f(x_k),y-x_k\rangle\right\}
    \\
    x_{k+1} &= (1-\tau_{k+1}) y_{k+1} + \tau_{k+1} z_{k+1}
\end{align*}
for $k=0,1,\dots$, where $x_0=z_0$ and $\{\alpha_k\}_{k=1}^\infty$ and $\{\tau_k\}_{k=1}^\infty$ are positive sequences to be determined.

We obtain AGM by performing a non-tight analysis of \eqref{eq:linear-coupling} and letting the analysis inform the choices of $\{\alpha_k\}_{k=1}^\infty$ and $\{\tau_k\}_{k=1}^\infty$.
The first step of this analysis is
\begin{align*}
\alpha_{k+1} \langle \nabla f(x_{k}), z_k - x_\star \rangle &\leq \frac{\alpha_{k+1}^2}{2} \norm{\nabla f(x_{k})}_*^2 + V_{z_k}(x_\star) - V_{z_{k+1}}(x_\star)
\\&\leq \alpha_{k+1}^2 L (f(x_{k}) - f(y_{k+1}))  + V_{z_k}(x_\star) - V_{z_{k+1}}(x_\star).
\end{align*}
The second inequality follows from 
\[
f(x_{k}) - f(y_{k+1}) \geq \frac{1}{2L} \norm{\nabla f(x_{k})}_*^2   +\underline{\frac{1}{2L} \norm{\nabla f(y_{k+1})}_*^2},
\]
but the underscored term $\frac{1}{2L} \norm{\nabla f(y_{k+1})}_*^2$ is not used,
i.e., proof utilizes the weaker and non-tight inequality
\[
f(x_{k}) - f(y_{k+1}) \geq \frac{1}{2L} \norm{\nabla f(x_{k})}_*^2.
\]
The second step of this analysis is to choose $\tau_k = \frac{1}{\alpha_{k+1}L}$ to eliminate $f(x_k)$ and to show
\begin{align*}
   \alpha_{k+1}^2L \bigl(&f(y_{k+1}) - f_\star \bigr)+ V_{z_{k+1}}(x_\star) \leq \left(\alpha_{k+1}^2L - \alpha_{k+1}\right)\left(f(y_{k}) - f_\star\right) + V_{z_k}(x_\star).
\end{align*}
The inequality follows from
\[
f(x_{k}) - f_\star \leq \langle \nabla f(x_{k}), x_{k+1} - x_\star \rangle -\underline{ \frac{1}{2L} \norm{\nabla f(x_{k})}_*^2}
\]
and 
\begin{align*}
&\langle \nabla f(x_{k}), y_k - x_{k} \rangle \leq f(y_k) - f(x_{k})- \underline{\frac{1}{2L}\norm{\nabla f(y_{k}) - \nabla f(x_{k})}_*^2},
\end{align*}
but the underscored terms are not used. Finally, convergence is established through a telescoping sum argument as Appendix \ref{sec:telescoping}.


\subsection{Linear coupling analysis of OGM}
We now derive OGM through performing a tight analysis of \eqref{eq:linear-coupling} and letting the analysis inform the choices of $\{\alpha_k\}_{k=1}^\infty$ and $\{\tau_k\}_{k=0}^\infty$.

In the first step of our linear coupling analysis, we follow the same arguments but do not take the step utilizing the non-tight inequality.
\begin{lemma}
\label{lemma:lc-first}
Assume (\hyperlink{A1}{A1}) and (\hyperlink{A2}{A2}). The iterates \eqref{eq:linear-coupling} satisfy
\begin{align*}
&\alpha_{k+1} \langle \nabla f(x_{k}), z_k - x_\star \rangle  \leq \frac{\alpha_{k+1}^2}{2} \norm{\nabla f(x_{k})}_*^2 + V_{z_k}(x_\star) - V_{z_{k+1}}(x_\star)
\end{align*}
for $k=0,1,\dots$.
\end{lemma}
\begin{proof}
This is exactly the first part of Lemma~4.2 of \cite{allen2014linear}.
\qed\end{proof}

In the second step of our linear coupling analysis, we choose $\tau_k = \frac{2}{\alpha_{k+1}L}$ to allow for a telescoping sum argument and show the following lemma.
\begin{lemma}
\label{lem:ogm-coupling-lemma}
Assume (\hyperlink{A1}{A1}), (\hyperlink{A2}{A2}) and (\hyperlink{A3}{A3}). 
Let $0<\tau_k = \frac{2}{\alpha_{k+1} L}\leq1$ for $k=0,1,..$, $\alpha_1 = \frac{2}{L}$, and $x_{-1} = x_0$.
Set 
$h(x) = f(x) - f_\star - \frac{1}{2L}\norm{\nabla f(x)}_*^2$.
The iterates \eqref{eq:linear-coupling} satisfy
\begin{align*}
    &\frac{\alpha_{k+1}^2L}{2} h(x_k) + V_{z_{k+1}}(x_\star)  \leq \frac{\alpha_{k+1}^2L - 2\alpha_{k+1}}{2} h(x_{k-1}) + V_{z_k}(x_\star)
\end{align*}
for $k=0,1,\dots$.
\end{lemma}
\begin{proof}

For $k= 1, 2, \dots$, we have
\begingroup
\allowdisplaybreaks
\begin{align}
    &\alpha_{k+1}\left(f(x_{k}) - f_\star)\right) \nonumber
    \\
    &\leq \alpha_{k+1}\langle \nabla f(x_{k}), x_{k} - x_\star \rangle - \frac{\alpha_{k+1}}{2L} \norm{\nabla f(x_{k})}_*^2 \label{lemma2-eq1}
    \\
    &= \alpha_{k+1}\langle \nabla f(x_{k}), x_{k} - z_k \rangle + \alpha_{k+1}\langle \nabla f(x_{k}), z_k - x_\star \rangle - \frac{\alpha_{k+1}}{2L} \norm{\nabla f(x_{k})}_*^2 \nonumber
    \\
    &= \frac{1-\tau_k}{\tau_k}\alpha_{k+1}\langle \nabla f(x_{k}), y_k - x_{k} \rangle + \alpha_{k+1}\langle \nabla f(x_{k}), z_k - x_\star \rangle - \frac{\alpha_{k+1}}{2L} \norm{\nabla f(x_{k})}_*^2 \nonumber
    \\
    &=  \frac{1-\tau_k}{\tau_k}\alpha_{k+1}\langle \nabla f(x_{k}), x_{k-1} - x_{k} - \frac{1}{L}Q^{-1}\nabla f(x_{k-1}) \rangle \nonumber
    \\
    &\qquad + \alpha_{k+1}\langle \nabla f(x_{k}), z_k - x_\star \rangle - \frac{\alpha_{k+1}}{2L} \norm{\nabla f(x_{k})}_*^2 \label{lemma2-eq2}
    \\
    &\leq  \frac{1-\tau_k}{\tau_k}\alpha_{k+1}\left( f(x_{k-1}) - f(x_{k}) - \frac{1}{2L} \norm{\nabla f(x_{k-1})}_*^2 - \frac{1}{2L}\norm{\nabla f(x_{k})}_*^2\right) \label{lemma2-eq3}
    \\
    &\qquad \qquad \qquad \qquad \qquad \qquad \qquad+ \alpha_{k+1}\langle \nabla f(x_{k}), z_k - x_\star \rangle - \frac{\alpha_{k+1}}{2L} \norm{\nabla f(x_{k})}_*^2 \nonumber 
    \\
    &\leq  \frac{1-\tau_k}{\tau_k}\alpha_{k+1}\left( f(x_{k-1}) - f(x_{k}) - \frac{1}{2L} \norm{\nabla f(x_{k-1})}_*^2 - \frac{1}{2L}\norm{\nabla f(x_{k})}_*^2\right) \label{eq:lemma2-final}
    \\
    &\qquad \qquad \qquad  +\frac{\alpha_{k+1}^2}{2} \norm{\nabla f(x_{k})}_*^2 + V_{z_k}(x_\star) - V_{z_{k+1}}(x_\star) - \frac{\alpha_{k+1}}{2L} \norm{\nabla f(x_{k})}_*^2. \nonumber
\end{align}
\endgroup
\eqref{lemma2-eq1} and \eqref{lemma2-eq3} follow from Lemma \ref{lem:cocoineq_generalnorm}, \eqref{lemma2-eq2} follows from the definition of linear coupling, and \eqref{eq:lemma2-final} follows from Lemma \ref{lemma:lc-first}.

The case of $k=0$ follows from $\alpha_1 = \frac{2}{L}$ and $f_\star -f(x_0) - \langle \nabla f(x_0), x_\star - x_0 \rangle - \frac{1}{2L} \norm{\nabla f(x_0)}_*^2 \geq 0$ with Lemma \ref{lemma:lc-first}.
\qed\end{proof}

\begin{theorem}
\label{thm:ogm-coupling-primary}
Assume (\hyperlink{A1}{A1}), (\hyperlink{A2}{A2}), and (\hyperlink{A3}{A3}).
Let the positive sequence $\{\alpha_k \}_{k=1}^\infty$ satisfy $0 \leq \alpha_{k+1}^2L - 2\alpha_{k+1} \leq \alpha_k^2 L$ for $k=1,2\dots$ and $\alpha_1 = \frac{2}{L}$.
Let $\tau_k = \frac{2}{\alpha_{k+1}L}$ for $k=1,2,\dots$.
The $y_k$-sequence of \eqref{eq:linear-coupling} exhibits the rate
\begin{align*}
    f(y_{k}) - f_\star \leq \frac{2V_{x_0}(x_\star)}{L\alpha_{k}^2}
\end{align*}
for $k=1,2, \dots$.
\end{theorem}
\begin{proof}
Sum the inequality of Lemma~\ref{lem:ogm-coupling-lemma} from $0$ to $(k-1)$.
Then use $V_{z_k}(x_\star) \geq 0$ and $f(y_{k}) \leq f(x_{k-1}) - \frac{1}{2L}\norm{\nabla f(x_{k-1})}_*^2$ to conclude the rate.
\qed\end{proof}

The $\{\theta_k\}_{k=0}^\infty$ of the original OGM formulation is related to  $\{\alpha_k\}_{k=1}^\infty$ through $\alpha_{k+1} = {2\theta_{k}}/{L}$ for $k =0,1,\dots$.
The seemingly different parameter choices $\tau_k=\frac{1}{\alpha_{k+1}L}$ for AGM and $\tau_k=\frac{2}{\alpha_{k+1}L}$ for OGM actually turn out to be the same as $\{\alpha_k\}_{k=1}^\infty$ for AGM and OGM differ by a factor of $2$.

The parameters $\{\alpha_{k}\}_{k=1}^\infty$ and $\{\tau_{k}\}_{k=1}^\infty$ are chosen to make the telescoping sum argument work and to make it work tightly, as described in Section~\ref{sec:telescoping}.
Specifically, one starts with a general form
\begin{align*}
    M_{k}\biggl(f(x_{k}) - f_\star - &B_{k}\norm{\nabla f(x_{k})}_*^2\biggr)  + V_{z_{k+1}}(x_\star)
    \\
    &\leq N_{k-1}\left(f(x_{k-1}) - f_\star - B_{k-1} \norm{\nabla f(x_{k-1})}_*^2\right) + V_{z_k}(x_\star),
\end{align*}
where the scalar coefficients $M_k$, $N_{k-1}$, $B_k$, and $B_{k-1}$ are determined by \eqref{eq:lemma2-final}.
To make the telescoping sum argument work as Appendix \ref{sec:telescoping}, the $\{B_k\}_{k=0}^\infty$ must be independent of $k$. 
The only possibility for this is $B_k=B_{k-1}=\frac{1}{2L}$, and this is achieved when
$$-\frac{1}{2L} \left(\alpha_{k+1} + \frac{1-\tau_k}{\tau_k} \alpha_{k+1}\right) = -\frac{\alpha_{k+1}^2}{2} + \frac{1}{2L} \left(\alpha_{k+1} + \frac{1-\tau_k}{\tau_k} \alpha_{k+1}\right) $$
holds. Solving this equation leads to the choice $\tau_k = \frac{2}{L\alpha_{k+1}}$.
The requirement $\alpha_{k+1}^2L-2\alpha_{k+1}\le \alpha_k^2L$ is needed for the telescoping sum argument to work, and the choice
$\alpha_{k+1}^2L-2\alpha_{k+1}= \alpha_k^2L$ makes the argument tight.

\subsection{Secondary sequence analysis}
In the linear coupling context, the last-step modification can be expressed as
\begin{equation}
    \tilde{x}_{k} = (1-\tilde{\tau}_k)y_k +\tilde{\tau}_k z_k 
    \label{eq:lc-last-step}
\end{equation}
for $k=0,1,\dots$, where $\{\tilde{\tau}_k\}_{k=0}^\infty$ is a positive sequence to be determined.

\begin{lemma}
\label{lem:ogm-coupling-last-step}
Assume (\hyperlink{A1}{A1}), (\hyperlink{A2}{A2}) and (\hyperlink{A3}{A3}). Let $0 < \tilde{\tau}_{k} = \frac{1}{\tilde{\alpha}_{k+1}L} \leq 1$ for $k=0,1,\dots$, $\tilde{\alpha}_1 = \frac{1}{L}$, and $x_{-1} = x_0$. 
Then the $\tilde{x}_k$-sequence of \eqref{eq:lc-last-step},
the secondary sequence with last-step modification of \eqref{eq:linear-coupling}, satisfies
\begin{align*}
    &\tilde{\alpha}_{k+1}^2L \left(f(\tilde{x}_{k})- f_\star\right)  + V_{z_{k+1}}(x_\star)  \leq \left(\tilde{\alpha}_{k+1}^2 L - \tilde{\alpha}_{k+1}\right) h(x_{k-1}) + V_{z_k}(x_\star)
\end{align*}
for  $k=0,1,\dots$.
\end{lemma}
\begin{proof}
Proof is identical to that of Lemma~\ref{lem:ogm-coupling-lemma} with substituted $\tau_k$ by $\tilde{\tau}_k$.
\qed\end{proof}

\begin{theorem}
\label{thm:ogm-coupling-secondary}
In the setup of Theorem~\ref{thm:ogm-coupling-primary},
let $0 \leq \tilde{\alpha}_{k+1}^2L - \tilde{\alpha}_{k+1} \leq \frac{1}{2}\alpha_k^2 L$ and $\tilde{\alpha}_1 = \frac{1}{L}$.
Then the $\tilde{x}_k$-sequence, the secondary sequence with last-step modification, of the linear coupling method \eqref{eq:linear-coupling} exhibits the rate
\begin{align*}
    f(\tilde{x}_{k}) - f_\star \leq \frac{V_{x_0}(x_\star)}{L\tilde{\alpha}_{k+1}^2}
\end{align*}
for $k=0,1,\dots$
\end{theorem}
\begin{proof}
Sum the inequality of Lemma~\ref{lem:ogm-coupling-lemma} from $0$
to $(k-2)$ and the inequality of Lemma~\ref{lem:ogm-coupling-last-step} with $k-1$.
Then use $V_{z_k} (x_\star) \geq 0$ to conclude the rate.
\qed\end{proof}

\subsection{Comparison of the linear coupling analyses of AGM and OGM}
The linear coupling analysis of Allen-Zhu and Orecchia \cite{allen2014linear}, which derives AGM, relies on the following two key lemmas.

\begin{lemma} \emph{\textbf{\cite[Lemma 4.2]{allen2014linear}}}
In the linear coupling setup,
\begin{align*}
    \alpha_{k+1} \langle \nabla f(x_k), z_k -x_\star \rangle &\leq \frac{\alpha_{k+1}^2}{2}\norm{\nabla f(x_{k})}_*^2 + V_{z_k}(x_\star) - V_{z_{k+1}}(x_\star)
    \\
    &\leq \alpha_{k+1}^2 L \left(f(x_k) - f(y_{k+1})\right) + V_{z_k}(x_\star) - V_{z_{k+1}}(x_\star)
\end{align*}
for $k=0,1,\dots$.
\end{lemma}

\begin{lemma} \emph{\textbf{\cite[Lemma 4.3]{allen2014linear}}}
\emph{(Coupling Lemma)}
In the linear coupling setup,
$$\alpha_{k+1}^2 L \left( f(y_{k+1}) -f_\star\right) + V_{z_{k+1}}(x_\star)\leq (\alpha_{k+1}^2 L - \alpha_{k+1})\left(f(y_{k}) - f_\star \right) +V_{z_{k}}(x_\star).$$ 
for $k=0,1, \dots$.
\end{lemma}

As discussed, the proof of \cite[Lemma 4.2]{allen2014linear} uses of the non-tight inequality 
\[
f(x_{k}) - f(y_{k+1}) \geq \frac{1}{2L} \norm{\nabla f(x_{k})}_*^2, 
\]
and the proof of \cite[Lemma 4.3]{allen2014linear} follows steps similar to that of Lemma~\ref{lem:ogm-coupling-lemma}, but uses the non-tight inequalities
\[
f(x_{k}) - f_\star \leq \langle \nabla f(x_{k}), x_{k+1} - x_\star \rangle 
\]
and 
\begin{align*}
&\langle \nabla f(x_{k}), y_k - x_{k} \rangle \leq f(y_k) - f(x_{k}).
\end{align*}

In both linear coupling analyses, for OGM and AGM, the telescoping sum argument is made tight by choosing $\{\alpha_k\}_{k=1}^\infty$ and $\{\tau_k\}_{k=1}^\infty$ appropriately.
However, the analysis of Allen-Zhu and Orecchia \cite{allen2014linear} uses non-tight inequalities before the telescoping sum argument, while our analysis uses tight inequalities in all steps.

\subsection{Unification of AGM and OGM}
\label{ss:agm-ogm-unification}
If we choose $w(y) = \frac{1}{2t} \norm{y}^2$, so that $V_x(y) = \frac{1}{2t}\norm{x-y}^2$, and $ 0 < t \leq 1$, so that $w$ is 1-strongly convex, and substitute $\alpha_{k+1} = {2\theta_{k}}/{L}$,
\eqref{eq:linear-coupling} becomes
\begin{align*}
  y_{k+1} &= x_{k} - \frac{1}{L}\nabla{f(x_{k})}
\\
  z_{k+1} &= z_{k} - \frac{2t\theta_k}{L}\nabla{f(x_{k})} 
\\
  x_{k+1} &= \left(1-\frac{1}{\theta_{k+1}}\right)y_{k+1} + \frac{1}{\theta_{k+1}}z_{k+1}
\end{align*}
for $k=0,1,\dots$.
We also express this method with the momentum and correction terms and without the $z^k$-iterates in Lemma \ref{AGM-OGM unification}.
This method unifies AGM and OGM through the constant $t$; AGM and OGM respectively correspond to $t= (1/2)$ and  $t=1$.

\begin{corollary}
\label{cor:agm-ogm-rate}
Assume (\hyperlink{A1}{A1}), (\hyperlink{A2}{A2}) and (\hyperlink{A3}{A3}). 
Let $0<t\le 1$. Then 
\[
    f(y_{k}) - f_\star \leq \frac{L \norm{x_0 - x_\star}^2}{4t\theta_{k-1}^2}
    \]
    for $k = 1,2, \dots$
\end{corollary}
\begin{proof}
This follows from Theorem~\ref{thm:ogm-coupling-primary} with 
$\alpha_{k+1} = \frac{2\theta_{k}}{L}$.
\qed\end{proof}
The rates of Corollary~\ref{cor:agm-ogm-rate} at $t= \frac{1}{2}$ and  $t=1$ exactly match the previously discussed rates of AGM and OGM.

\begin{lemma}
\label{AGM-OGM unification}
The unified form is equivalent to
\begin{align*}
  y_{k+1} &= x_{k} - \frac{1}{L}\nabla{f(x_{k})}
\\
  x_{k+1} &= y_{k+1} + \frac{\theta_{k}-1}{\theta_{k+1}}(y_{k+1}- y_k)+(2t-1)\frac{\theta_k}{\theta_{k+1}}(y_{k+1} - x_{k}).
\end{align*}
\end{lemma}
\begin{proof}
To prove the equivalency, we show that the above sequence leads to $$x_{k+1} = \left(1- \frac{1}{\theta_{k+1}}\right)y_{k+1} + \frac{1}{\theta_{k+1}}z_{k+1}.$$
That is,
\begingroup
\allowdisplaybreaks
\begin{align*}
    x_{k+1} &= \left(1-\frac{1}{\theta_{k+1}}\right) y_{k+1} + \frac{\theta_k}{\theta_{k+1}}y_{k+1}-  \frac{\theta_k -1}{\theta_{k+1}}y_{k} - (2t-1)\frac{\theta_k}{\theta_{k+1}}\frac{1}{L}\nabla f(x_k)
    \\
    &= \left(1-\frac{1}{\theta_{k+1}}\right)y_{k+1} + \frac{\theta_k}{\theta_{k+1}}\left(x_{k}-\frac{1}{L}\nabla f(x_k)\right)
    \\ &\qquad -  \frac{\theta_k -1}{\theta_{k+1}}y_{k} - (2t-1)\frac{\theta_k}{\theta_{k+1}}\frac{1}{L}\nabla f(x_k)
    \\
    &= \left(1-\frac{1}{\theta_{k+1}}\right)y_{k+1} + \frac{\theta_k}{\theta_{k+1}}x_{k}-  \frac{\theta_k -1}{\theta_{k+1}}y_{k} - 2t\frac{\theta_k}{\theta_{k+1}}\frac{1}{L}\nabla f(x_k) 
    \\
    &= \left(1-\frac{1}{\theta_{k+1}}\right)y_{k+1} -  \frac{\theta_k -1}{\theta_{k+1}}y_{k} - 2t\frac{\theta_k}{\theta_{k+1}}\frac{1}{L}\nabla f(x_k) 
    \\
    &\qquad + \frac{\theta_k}{\theta_{k+1}}\left(y_{k} + \frac{\theta_{k-1} -1}{\theta_{k}}(y_{k}- y_{k-1}) - (2t-1) \frac{\theta_{k-1}}{\theta_{k}}\frac{1}{L}\nabla f(x_{k-1})\right)
    \\
    &= \left(1-\frac{1}{\theta_{k+1}}\right)y_{k+1} + \left(\frac{\theta_k}{\theta_{k+1}}  + \frac{\theta_{k-1}-1}{\theta_{k+1}}-\frac{\theta_{k}-1}{\theta_{k+1}}\right)y_{k} - \frac{\theta_{k-1} -1}{\theta_{k+1}}y_{k-1}
    \\
    &\qquad- (2t-1)\frac{\theta_{k-1}}{\theta_{k+1}}\frac{1}{L}\nabla f(x_{k-1}) - 2t\frac{\theta_k}{\theta_{k+1}}\frac{1}{L}\nabla f(x_k) 
    \\
    &= \left(1-\frac{1}{\theta_{k+1}}\right)y_{k+1} + \frac{\theta_{k-1}}{\theta_{k+1}}y_{k} - \frac{\theta_{k-1} -1}{\theta_{k+1}}y_{k-1} 
    \\
    &\qquad - (2t-1)\frac{\theta_{k-1}}{\theta_{k+1}}\frac{1}{L}\nabla f(x_{k-1}) - 2t\frac{\theta_k}{\theta_{k+1}}\frac{1}{L}\nabla f(x_k)
    \\
    &= \left(1-\frac{1}{\theta_{k+1}}\right)y_{k+1} + \frac{\theta_{k-1}}{\theta_{k+1}}\left(x_{k-1}- \frac{1}{L}\nabla f(x_{k-1})\right) - \frac{\theta_{k-1} -1}{\theta_{k+1}}y_{k-1} 
    \\
    &\qquad - (2t-1)\frac{\theta_{k-1}}{\theta_{k+1}}\frac{1}{L}\nabla f(x_{k-1}) - 2t\frac{\theta_k}{\theta_{k+1}}\frac{1}{L}\nabla f(x_k)    
    \\
    &= \left(1-\frac{1}{\theta_{k+1}}\right)y_{k+1} + \frac{\theta_{k-1}}{\theta_{k+1}}x_{k-1}- \frac{\theta_{k-1} -1}{\theta_{k+1}}y_{k-1} 
    \\&\qquad- 2t\frac{\theta_k}{\theta_{k+1}}\frac{1}{L}\nabla f(x_k) -2t \frac{\theta_{k-1}}{\theta_{k+1}}\frac{1}{L}\nabla f(x_{k-1})   
    \\
    & \qquad \qquad \qquad \qquad\qquad  \vdots 
    \\
    &= \left(1-\frac{1}{\theta_{k+1}}\right)y_{k+1} + \frac{\theta_{0}}{\theta_{k+1}}x_{0}- \frac{\theta_{0} -1}{\theta_{k+1}}y_{0} - \frac{1}{\theta_{k+1}}\sum_{i=0}^{k}{2t\theta_i \frac{1}{L}\nabla f(x_i)}
    \\
    &= \left(1- \frac{1}{\theta_{k+1}}\right)y_{k+1} + \frac{1}{\theta_{k+1}}z_{k+1}.
\end{align*}
\endgroup
\qed\end{proof}

\subsection{Discussion}
By identifying OGM as an instance of linear coupling, we generalized OGM to the setup with quadratic norms and mirror descent steps while maintaining the factor-$\sqrt{2}$ improvement.
However, we do point out that the generalization is narrower than that of \cite{allen2014linear}, which allows non-quadratic norms and constrained $y_k$-and $z_k$-updates.
The analysis on strongly convex case follows from a similar line of reasoning, and is presented in Appendix, Section \ref{sec:lc-sc-ogm}.


In addition to the human-understandable derivation of OGM, this section provides two non-obvious observations, which we point out again.
The first is that AGM and OGM can be unified into a single parameterized family of accelerated gradient methods, all achieving the $\mathcal{O}({1/k^2})$ rate.
Another is that the linear coupling analysis of Allen-Zhu and Orecchia \cite{allen2014linear} was suboptimal in the same way that AGM is suboptimal and can be improved.

\section{Conclusion}
In this work, we presented human-understandable analyses of OGM.
The first key insight is to use a Lyapunov function with $f(x_k) - f_\star - \frac{1}{2L}\norm{\nabla{f(x_k)}}^2$, a somewhat unusual term in Lyapunov analyses.
The second key insight is to obtain OGM by refining the linear coupling analysis of Allen-Zhu and Orecchia \cite{allen2014linear} through replacing non-tight inequalities with tight ones.
With these insights, we extended the factor-$\sqrt{2}$ acceleration to other setups.

In our view, the most significant contribution of this work is the improved understanding of OGM's acceleration mechanism.
While Nesterov's acceleration mechanism has been utilized as a component in a wide range of setups, OGM's acceleration mechanism has not yet seen any external use.
Through the understanding provided by the analysis of this work, we hope OGM's acceleration becomes more widely utilized to gain a (perhaps factor-$\sqrt{2}$) speedup compared to what can be achieved with AGM's acceleration.
For example, whether accelerated coordinate gradient methods \cite{allen2016even, nesterov2017efficiency} or non-convex stochastic optimization \cite{ghadimi2016accelerated} can be improved with OGM's acceleration mechanism would be an interesting question to address in future work.
Improving the FISTA \cite{beck2009fast} and the more general mirror descent setup \cite{doi:10.1287/moor.2016.0817,doi:10.1137/16M1099546} are also interesting directions, although there are known limitations \cite{kim2018fista,Dragomir2021}.


Finally, studying how OGM's acceleration interacts with other techniques used to analyze AGM, such as the continuous-time analysis \cite{su2014differential}, high-resolution ODEs \cite{shi2019acceleration}, and variational perspective \cite{WibisonoE7351} is also an interesting direction.

\section*{Acknowledgements}
CP was supported by an undergraduate research internship in the second half of the 2020 Seoul National University College of Natural Sciences. 
JP and EKR were supported by the National Research Foundation of Korea (NRF) Grant funded by the Korean Government (MSIP) [No. 2020R1F1A1A01072877], the National Research Foundation of Korea (NRF) Grant funded by the Korean Government (MSIP) [No. 2017R1A5A1015626], and by the New Faculty Startup Fund from Seoul National University.
We thank Gyumin Roh for reviewing the manuscript and providing valuable feedback. 
We thank Bryan Van Scoy and Suvrit Sra for the discussions regarding the triple momentum method and estimate sequences, respectively.

\clearpage


%
\section*{Conflict of interest}
 The authors declare that they have no conflict of interest.

\bibliographystyle{spmpsci}      
\bibliography{faagm}
\clearpage
\appendix

\section{Method reference}
For reference, we restate all aforementioned methods. In all methods, we assume that $f$ is $L$-smooth function, $\{ \theta_k \}_{k=0}^\infty$ and $\{ \varphi_k \}_{k=0}^\infty$ are the sequences of positive scalars, and $x_0=y_0=z_0$. 

\paragraph{OGM.} One form of OGM is  
\begin{align*}
  y_{k+1} &= x_{k} - \frac{1}{L}\nabla{f(x_{k})}
\\
  x_{k+1} &= y_{k+1} + \frac{\theta_{k}-1}{\theta_{k+1}}(y_{k+1}- y_k) + \frac{\theta_k}{\theta_{k+1}}(y_{k+1} - x_{k})
\end{align*}
and an equivalent form with $z$-iterates is
\begin{align*}
  y_{k+1} &= x_{k} - \frac{1}{L}\nabla{f(x_{k})}
\\
  z_{k+1} &= z_{k} - \frac{2\theta_k}{L}\nabla{f(x_{k})} 
\\
  x_{k+1} &= \left(1-\frac{1}{\theta_{k+1}}\right)y_{k+1} + \frac{1}{\theta_{k+1}}z_{k+1}
\end{align*}
for $k=0,1,\dots$.
The \emph{last-step modification} on the secondary sequence can be written as 
\begin{align*}
\tilde{x}_{k+1} &= y_{k+1} + \frac{\theta_{k}-1}{\varphi_{k+1}}(y_{k+1}- y_k) + \frac{\theta_k}{\varphi_{k+1}}(y_{k+1} - x_{k})\\
    &= \left(1-\frac{1}{\varphi_{k+1}}\right)y_{k+1} + \frac{1}{\varphi_{k+1}}z_{k+1}
\end{align*}
where $k=0,1,\dots$.
\medskip

\paragraph{OGM-simple.} OGM-simple is a simpler variant of \textbf{OGM} with $\theta_k = \frac{k+2}{2}$ and $\varphi_k = \frac{k+1 +\frac{1}{\sqrt{2}}}{\sqrt{2}}$. One form of OGM-simple is
\begin{align*}
    y_{k+1} &= x_{k} - \frac{1}{L}\nabla{f(x_{k})}
\\
  x_{k+1} &= y_{k+1} + \frac{k}{k+3}(y_{k+1}- y_k) + \frac{k+2}{k+3}(y_{k+1} - x_{k})
\end{align*}
and an equivalent form with $z$-iterates is
\begin{align*}
  y_{k+1} &= x_{k} - \frac{1}{L}\nabla{f(x_{k})}
\\
  z_{k+1} &= z_{k} - \frac{k+2}{L}\nabla{f(x_{k})} 
\\
  x_{k+1} &= \left( 1-\frac{2}{{k+3}} \right)y_{k+1} + \frac{2}{k+3}z_{k+1}
\end{align*}
for $k=0,1,\dots$. The \emph{last-step modification} on secondary sequence is written as
\begin{align*}
  \tilde{x}_{k+1}&= y_{k+1} + \frac{k}{\sqrt{2}(k+2) + 1}(y_{k+1} - y_{k}) + \frac{k+2}{\sqrt{2}(k+2) + 1}(y_{k+1} - x_k)
\end{align*}
where $k=0,1,\dots$.
\medskip

\paragraph{SC-OGM.}
Here, we assume that $f$ is a $\mu$-strongly convex function, condition number of $f$ is $\kappa = L/\mu$, and $\gamma = \frac{\sqrt{8\kappa+1}+3}{2\kappa -2}$.
SC-OGM is written as
\begin{align*}
    y_{k+1} &= x_{k} - \frac{1}{L}\nabla{f(x_{k})}
\\
    x_{k+1} &= y_{k+1} + \frac{1}{2\gamma + 1}(y_{k+1}-y_k) + \frac{1}{2\gamma + 1}(y_{k+1} - x_{k})
\end{align*}
for $k=0,1,\dots$.
\medskip

\paragraph{LC-OGM.} 
LC-OGM (Linear Coupling OGM) is defined as
\begin{align*}
    y_{k+1} &= x_k-L^{-1}Q^{-1}\nabla f(x_k)
    \\
    z_{k+1} &= \argmin_{y\in\mathbb{R}^n}\left\{V_{z_k}(y)+\langle \alpha_{k+1}\nabla f(x_k),y-x_k\rangle\right\}
    \\
    x_{k+1} &= (1-\tau_{k+1}) y_{k+1} + \tau_{k+1} z_{k+1}
\end{align*}
for $k=0,1,\dots$, where $V_z(y)$ is a Bregman divergence, $\{\alpha_k\}_{k=1}^\infty$ and $\{\tau_k\}_{k=1}^\infty$ are nonnegative sequences defined as $\alpha_1 = \frac{2}{L}$, $0 \leq \alpha_{k+1}^2L -2\alpha_{k+1} \leq \alpha_k^2L$, $\tau_k = \frac{2}{\alpha_{k+1}L}$, and $Q$ is a positive definite matrix defining $\norm{x}^2 = x^T Q x $.

For \emph{last step modification}, we define positive sequences $\{\tilde{\alpha}_k\}_{k=1}^\infty$ and $\{\tilde{\tau}_k\}_{k=1}^\infty$ as $\alpha_1 = \frac{1}{L}$, $0 \leq \tilde{\alpha}_{k+1}^2 L - \tilde{\alpha_{k+1}} \leq \frac{1}{2}\alpha_k^2L$, and $\tilde{\tau}_k = \frac{1}{\tilde{\alpha}_{k+1}L} $, and also define 
\[
\tilde{x}_k = (1-\tilde{\tau}_k)y_k + \tilde{\tau}_k z_k
\]
for $k=1,2,\dots$.
\medskip

\paragraph{Unification of AGM and OGM.}
Using LC-OGM, we can unify AGM and OGM as
\begin{align*}
  y_{k+1} &= x_{k} - \frac{1}{L}\nabla{f(x_{k})}
\\
  z_{k+1} &= z_{k} - \frac{2t\theta_k}{L}\nabla{f(x_{k})} 
\\
  x_{k+1} &= \left(1-\frac{1}{\theta_{k+1}}\right)y_{k+1} + \frac{1}{\theta_{k+1}}z_{k+1}.
\end{align*}
for $k=0,1,\dots$.
This is equivalent to 
\begin{align*}
  y_{k+1} &= x_{k} - \frac{1}{L}\nabla{f(x_{k})}
\\
  x_{k+1} &= y_{k+1} + \frac{\theta_{k}-1}{\theta_{k+1}}(y_{k+1}- y_k)+(2t-1)\frac{\theta_k}{\theta_{k+1}}(y_{k+1} - x_{k}).
\end{align*}
\medskip

\paragraph{LC-SC-OGM.} LC-SC-OGM (Linear Coupling Strongly Convex OGM) is
\begin{align*}
    y_{k+1} &= x_{k} - \frac{1}{L}Q^{-1}\nabla{f(x_{k})}
\\
z_{k+1} &= \frac{1}{1+\gamma}\left(z_k + \gamma x_{k} -\frac{\gamma}{\mu}Q^{-1}\nabla f(x_{k})\right)
\\
    x_{k+1} &= \tau z_{k+1} + (1-\tau) y_{k+1},
\end{align*}
for $k=0,1,\dots$, where $Q$ is a positive definite matrix. 
\medskip

\section{Co-coercivity inequality in general norm}
\begin{lemma}
\label{lem:ccpfunction_lemma3}
Let $f$ be a closed convex proper function. Then, 
$$0 \leq f(x) + f^* (u) - \langle x, u \rangle$$ and 
$$\inf_{x} \{f(x) + f^*(u) - \langle x, u \rangle\} = 0 $$
$$\inf_{u} \{f(x) + f^*(u) - \langle x, u \rangle \} = 0.$$
\end{lemma}
\begin{proof}
By the definition of the conjugate function,
$$-f^*(u) = \inf_x \left \{f(x) - \langle x, u \rangle\right \} $$ 
and
$$\inf_{x} \{f(x) + f^*(u) - \langle x, u \rangle\} = 0. $$
Therefore,
$$0 \leq f(x) + f^* (u) - \langle x, u \rangle \quad \forall x .$$
The statement with $u$ follows from the same argument and the fact that $f^{**} = f$.
\qed\end{proof}

\begin{lemma} 
\label{lem:ineq_norm}
Consider a norm $\|\cdot\|$ and its dual norm $\|\cdot\|_*$.
Then,
$$0 \leq \frac{1}{2}\norm{x}^2 + \frac{1}{2} \norm{u}_*^2 - \langle x, u \rangle $$ 
and 
$$\inf_{x \in \mathbb{R}^n} \left \{\frac{1}{2}\norm{x}^2 + \frac{1}{2} \norm{u}_*^2 - \langle x, u \rangle \right \} =0$$ 
$$\inf_{u \in \mathbb{R}^n} \left \{\frac{1}{2}\norm{x}^2 + \frac{1}{2} \norm{u}_*^2 - \langle x, u \rangle \right \} = 0.$$
\end{lemma}
\begin{proof}
This follows from Lemma \ref{lem:ccpfunction_lemma3}
with $f(x) = \frac{1}{2}\norm{x}^2$ and $\left(\frac{1}{2}\norm{\cdot}^2\right)^* = \frac{1}{2}\norm{\cdot}_*^2$.
\qed\end{proof}

\begin{lemma}
\label{lem:grad}
Let 
$$\mathtt{Grad}(x) = \argmin_{y \in \mathbb{R}^n}\left \{ \frac{L}{2}\norm{y-x}^2 + \langle \nabla f(x), y-x \rangle \right \}.$$ Then, 
$$\langle \nabla f(x) , \mathtt{Grad}(x) - x \rangle + \frac{L}{2}\norm{\mathtt{Grad}(x) - x }^2 = -\frac{1}{2L}\norm{\nabla f(x)}_*^2.$$
\end{lemma}

\begin{proof}
Let $z= L (\mathtt{Grad}(x) - x)$.
By the definition of $\mathtt{Grad}(x)$ and Lemma \ref{lem:ineq_norm}, we have
\begin{align*}
    \frac{1}{2L}\norm{\nabla f(x)}_*^2 + &\frac{L}{2} \norm{\mathtt{Grad}(x) - x}^2 + \langle \nabla f(x) , \mathtt{Grad}(x) - x \rangle 
    \\
    &= \inf_{z \in \mathbb{R}^n} \frac{1}{2L}\norm{\nabla f(x)}_*^2 + \frac{1}{2L} \norm{z}^2 + \frac{1}{L} \langle \nabla f(x), z \rangle \\
&= 0.
\end{align*}
\qed\end{proof}

\begin{lemma}
Let $f: \mathbb{R}^n \rightarrow \mathbb{R}$ be a differentiable convex function such that $$\norm{\nabla f(x) - \nabla f(y)}_* \leq L \norm{x-y}$$
for all $x,y\in \mathbb{R}^n$.
Then
$$f(y) \leq f(x) + \langle \nabla f(x), y-x \rangle + \frac{L}{2}\norm{y-x}^2.$$
\end{lemma}

\begin{proof}
Since a differentiable convex function is continuously differentiable \cite[Theorem~25.5]{rockafellar1970},
\begin{align*}
f(y) - f(x) &= \int_0^1 \langle \nabla f(x + t(y-x)) , y-x \rangle dt
\\
&=  \int_0^1 \langle \nabla f(x + t(y-x)) - \nabla f(x) , y-x \rangle dt + \langle \nabla f(x), y-x \rangle 
\\
&\leq \int_0^1 \norm{\nabla f(x + t(y-x)) - \nabla f(x)}_* \norm{y-x} dt + \langle \nabla f(x), y-x \rangle
\\
&\leq \int_0^1 t L \norm{y-x}^2 dt + \langle \nabla f(x), y-x \rangle = \frac{L}{2} \norm{y-x}^2  + \langle \nabla f(x) , y-x \rangle .
\end{align*}
\qed\end{proof}

\begin{lemma}[Co-coercivity inequality with general norm]
\label{lem:cocoineq_generalnorm}
Let $f: \mathbb{R}^n \rightarrow \mathbb{R}$ be a differentiable convex function such that $$\norm{\nabla f(x) - \nabla f(y)}_* \leq L \norm{x-y}$$
for all $x,y\in \mathbb{R}^n$.
Then
$$f(y) \geq f(x) + \langle \nabla f(x), y-x \rangle + \frac{1}{2L}\norm{\nabla f(x)- \nabla f(y)}_*^2.$$
\end{lemma}

\begin{proof}
Set $\phi(y) = f(y) - \langle \nabla f(x), y-x\rangle $. Then $x \in \argmin \phi$. So by Lemma \ref{lem:grad},
\begin{align*}
    \phi(x) &\leq \phi(\mathtt{Grad}(y))
    \\
    &\leq \phi(y) + \langle \nabla \phi(y), \mathtt{Grad}(y) -y \rangle + \frac{L}{2}\norm{\mathtt{Grad}(y) - y}^2 
    \\
    &= \phi(y) - \frac{1}{2L}\norm{\nabla \phi(y)}_*^2. 
\end{align*}
Substituting $f$ back in $\phi$ yields the co-coercivity inequality. 
\qed\end{proof}

\section{Telescoping sum argument}
\label{sec:telescoping}
Suppose we established the inequality
\[
    a_i F_i + b_i G_i \leq c_i F_{i-1} + d_i G_{i-1} - E_i 
\]
for $i=1,2,\dots$, where $E_i, F_i, G_i$ are nonnegative quantities and $a_i$, $b_i$, $c_i$, and $d_i$ are nonnegative scalars.
Assume $c_i \leq a_{i-1}$ and $d_i \leq b_{i-1}$.
By summing the inequalities for $i=1, 2, \dots, k$, we obtain
\begin{align*}
    a_k F_k
    &\leq - b_k G_k
    - \sum_{i=2}^{k} (a_{i-1} - c_i) F_{i-1} - \sum_{i=2}^k (b_{i-1} - d_i) G_{i-1} - \sum_{i=2}^k E_i +    c_1 F_{0} + d_1 G_{0}\\
    &\le   c_1 F_{0} + d_1 G_{0}.
\end{align*}

However, note that the 
\[
- b_k G_k
    - \sum_{i=2}^{k} (a_{i-1} - c_i) F_{i-1} - \sum_{i=2}^k (b_{i-1} - d_i) G_{i-1} - \sum_{i=1}^k E_i
\]
terms are wasted in the analysis. If one has the freedom to do so, it may be good to choose parameters so that
\[
    a_{i-1} = c_i,\, b_{i-1} = d_i
\]
and $E_i = 0$ for $i=1,2,\dots$. Not having wasted terms may be an indication that the analysis is tight.


\section{SC-OGM via linear coupling}
\label{sec:lc-sc-ogm}
In this section, we analyze SC-OGM through the linear coupling analysis.
We consider the linear coupling form
\begin{align*}
    y_{k+1} &= x_{k} - \frac{1}{L}Q^{-1}\nabla{f(x_{k})}
\\
z_{k+1} &= \frac{1}{1+\gamma}\left(z_k + \gamma x_{k} -\frac{\gamma}{\mu}Q^{-1}\nabla f(x_{k})\right)
\\
    x_{k+1} &= \tau z_{k+1} + (1-\tau) y_{k+1},
\end{align*}
where $\tau$ is a coupling coefficient to be determined.
As an aside, we can view $z_{k+1}$ as a mirror descent update of the form
$$z_{k+1} =\argmin_{z} \left\{ \frac{1}{2}\norm{z-z_k}^2 + \frac{\gamma}{2}\norm{z-x_k}^2 +\frac{\gamma}{\mu} \langle \nabla f(x_k), z \rangle \right \}, $$
which is similar to what was considered in \cite{allen2016even}.

\begin{lemma}
\label{lem:lem5-4-analogue}
 Assume (\hyperlink{A1}{A1}), (\hyperlink{A2}{A2}) and (\hyperlink{A3}{A3}). Then, 
 \begin{align*}
 \frac{\gamma}{\mu}\langle \nabla f(x_{k}) , z_{k+1}& - x_\star \rangle - \frac{\gamma}{2}\norm{x_k - x_\star}^2 
 \\
 &\leq  -\frac{\gamma^2}{2(1+\gamma)\mu^2}\norm{\nabla f(x_k)}_*^2+ \frac{1}{2} \norm{z_k - x_\star}^2 - \frac{1+\gamma}{2} \norm{z_{k+1} - x_\star}^2 
 \end{align*}
for $k = 0,1, \dots$.
\end{lemma}
\begin{proof}
This proof follows steps similar to that of \cite[Lemma 5.4]{allen2016even}.

From the definition of $z_{k+1}$, we say 
\begin{align*}
0=& \langle \frac{\partial}{\partial z}\left\{ \frac{1}{2}\norm{z-z_k}^2 + \frac{\gamma}{2}\norm{z-x_k}^2 +\frac{\gamma}{\mu} \langle \nabla f(x_k), z \rangle \right \}\bigg|_{z_{k+1}}, z_{k+1} - x_\star \rangle
\\
= &\langle Q(z_{k+1} - z_k), z_{k+1} - x_\star \rangle + \frac{\gamma}{\mu} \langle \nabla f(x_{k}), z_{k+1} - x_\star \rangle + \gamma \langle Q(z_{k+1} - x_k),z_{k+1} - x_\star \rangle
\end{align*}
By three point equation, 
\begin{align*}
\frac{\gamma}{\mu}\langle \nabla f(x_{k}) ,z_{k+1}& - x_\star \rangle + \gamma \left( \frac{1}{2} \norm{x_k - z_{k+1}}^2 - \frac{1}{2} \norm{x_k - x_\star}^2 \right) 
\\
&= -\frac{1}{2} \norm{z_{k}- z_{k+1}}^2 + \frac{1}{2} \norm{z_k - x_\star}^2 - \frac{1+\gamma}{2} \norm{z_{k+1} - x_\star}^2. 
\end{align*}
Plugging the definition of $z_{k+1}$, 
\begin{align*}
\frac{\gamma}{2}& \norm{x_k - z_{k+1}}^2 + \frac{1}{2}\norm{z_k - z_{k+1}}^2 \\&= \frac{\gamma}{2} \norm{\frac{1}{1+\gamma}(x_k - z_k) + \frac{\gamma}{(1+\gamma)\mu} Q^{-1}\nabla f(x_k)}^2 + \frac{1}{2}\norm{-\frac{\gamma}{1+\gamma}(x_k - z_k) + \frac{\gamma}{(1+\gamma)\mu}Q^{-1}\nabla f(x_k)}^2
\\
&\geq \frac{\gamma^2}{2(1+\gamma)\mu^2}\norm{\nabla f(x_k)}_*^2.
\end{align*}
Combining results above, we get
 \begin{align*}
 \frac{\gamma}{\mu}\langle \nabla f(x_{k}) , z_{k+1}& - x_\star \rangle - \frac{\gamma}{2}\norm{x_k - x_\star}^2 
 \\
 &\leq  -\frac{\gamma^2}{2(1+\gamma)\mu^2}\norm{\nabla f(x_k)}_*^2+ \frac{1}{2} \norm{z_k - x_\star}^2 - \frac{1+\gamma}{2} \norm{z_{k+1} - x_\star}^2 .
 \end{align*}
\qed\end{proof}

\begin{lemma}[Coupling lemma in SC-OGM]
 Assume (\hyperlink{A1}{A1}), (\hyperlink{A2}{A2}) and (\hyperlink{A3}{A3}). Then 
 \begin{align*}
 (1+\gamma)\biggl(f(x_k) - \frac{1}{2L}&\norm{\nabla f(x_k)}_*^2 + \frac{\mu}{2}\norm{z_k - x_\star}^2 \biggr) \\&\leq \left(f(x_{k-1}) - \frac{1}{2L}\norm{\nabla f(x_{k-1})}_*^2 + \frac{\mu}{2}\norm{z_{k-1} - x_\star}^2 \right)
 \end{align*}
 holds for $k = 1, 2, \dots$
\end{lemma}
\begin{proof}
We have
\begingroup
\allowdisplaybreaks
\begin{align*}
    \gamma &\left( f(x_{k}) - f(x_\star) \right)
    \\
    &\leq \gamma \langle \nabla f(x_{k}), x_{k}- x_\star\rangle - \frac{\mu \gamma}{2}\norm{x_k - x_\star}^2 
    \\
    &= \gamma \langle \nabla f(x_{k}), x_{k}- z_{k}\rangle + \gamma \langle \nabla f(x_{k}), z_{k}- x_\star\rangle  - \frac{\mu \gamma}{2}\norm{x_k - x_\star}^2 
    \\
    &= \frac{1-\tau}{\tau}\gamma \langle \nabla f(x_k) ,y_k - x_k \rangle +  \gamma \langle \nabla f(x_{k}), z_{k}- x_\star\rangle  - \frac{\mu \gamma}{2}\norm{x_k - x_\star}^2 
    \\
    &= \frac{1-\tau}{\tau}\gamma \langle \nabla f(x_k) , x_{k-1} - x_{k} - \frac{1}{L}Q^{-1}\nabla f(x_{k-1}) \rangle +  \gamma \langle \nabla f(x_{k}), z_{k}- x_\star\rangle  - \frac{\mu \gamma}{2}\norm{x_k - x_\star}^2 
    \\
    & \leq \left(\frac{1-\tau}{\tau}\gamma - 1 \right)\langle \nabla f(x_k) , x_{k-1} - x_{k} - \frac{1}{L}Q^{-1}\nabla f(x_{k-1}) \rangle 
    \\
    &\qquad+\left( f(x_{k-1}) - f(x_k) - \frac{1}{2L}\norm{\nabla f(x_{k-1})}_*^2 - \frac{1}{2L}\norm{\nabla f(x_k)}_*^2\right)
    \\&\qquad \qquad +\gamma \langle \nabla f(x_{k}), z_{k}-z_{k+1} \rangle  + \gamma \langle \nabla f(x_{k}), z_{k+1}- x_\star\rangle   - \frac{\mu \gamma}{2}\norm{x_k - x_\star}^2 
    \\
    & \leq \left(\frac{1-\tau}{\tau}\gamma - 1 \right)\langle \nabla f(x_k) , y_{k} - x_{k}\rangle +\left( f(x_{k-1}) - f(x_k) - \frac{1}{2L}\norm{\nabla f(x_{k-1})}_*^2 - \frac{1}{2L}\norm{\nabla f(x_k)}_*^2\right)
    \\&\qquad +\gamma \langle \nabla f(x_{k}), z_{k}-z_{k+1} \rangle -\frac{\gamma^2}{2(1+\gamma)\mu}\norm{\nabla f(x_k)}_*^2+ \frac{\mu}{2} \norm{z_k - x_\star}^2 - \frac{(1+\gamma)\mu}{2} \norm{z_{k+1} - x_\star}^2,
\end{align*}
\endgroup
where the last inequality is an application of Lemma \ref{lem:lem5-4-analogue}.
Note that
\begin{align*}
    z_k - z_{k+1} &= z_k - \frac{1}{1+\gamma}\left(z_k + \gamma  x_k -\frac{\gamma}{\mu}Q^{-1}\nabla f(x_k) \right)
    \\
    &= \frac{\gamma}{1+\gamma}(z_k-x_k) + \frac{\gamma}{(1+\gamma)\mu} Q^{-1} \nabla f(x_k) 
    \\
    &= \frac{\gamma}{1+\gamma}\frac{1-\tau}{\tau}(x_k - y_k) + \frac{\gamma}{(1+\gamma)\mu} Q^{-1} \nabla f(x_k).
\end{align*}

To eliminate the $\langle \nabla f(x_k), \cdot \rangle$ term, we choose $\tau$ to satisfy
\begin{align} \frac{1-\tau}{\tau}\gamma - 1  = \frac{\gamma}{1+\gamma}\frac{1-\tau}{\tau}. \label{eq:lcsc-eq1}
\end{align}
Plugging this in, the inequality above is
\begin{align*}
    \gamma &\left( f(x_{k}) - f(x_\star) \right)
    \\
    &\leq 
    \left( f(x_{k-1}) - f(x_k) - \frac{1}{2L}\norm{\nabla f(x_{k-1})}_*^2 - \frac{1}{2L}\norm{\nabla f(x_k)}_*^2\right) 
    \\
    &\quad\quad 
    +\frac{\gamma^2}{2(1+\gamma)\mu}\norm{\nabla f(x_k)}_*^2+ \frac{\mu}{2} \norm{z_k - x_\star}^2 - \frac{(1+\gamma)\mu}{2} \norm{z_{k+1} - x_\star}^2.
\end{align*}
In order to make the telescoping form such as 
\begin{align*}
    M_{k}\biggl(f(x_{k})  - B_{k}\norm{\nabla f(x_{k})}_*^2  +& C_k \norm{z_{k+1} - x_\star}^2 \biggr)
    \\ &\leq N_{k-1}\left(f(x_{k-1})  - B_{k-1} \norm{\nabla f(x_{k-1})}_*^2+ C_{k-1} \norm{z_{k} - x_\star}^2 \right),
\end{align*}
we chose $B_k = \frac{1}{2L}$ and $C_k = \frac{\mu}{2}$, which leads to the choice of $\gamma$ satisfying
\begin{align}
    \frac{2+\gamma}{2L} = \frac{\gamma^2}{2(1+\gamma)\mu}. \label{eq:lcsc-eq2}
\end{align}
We get the desired result by plugging \eqref{eq:lcsc-eq1} and \eqref{eq:lcsc-eq2} in the above inequality.
\qed\end{proof}

\section{Asymptotic characterization of $\theta_k$}
\begin{theorem}
\label{thm:theta-asymp}
Let the positive sequence $\{\theta_k\}_{k=0}^\infty$ satisfy $\theta_0 = 1$ and $\theta_{k+1}^2 - \theta_{k+1} - \theta_{k}^2=0$ for $k = 0,1, \dots$. 
Then, 
\[
\theta_k = \frac{k+\zeta+1}{2} +  \frac{\log k}{4} + o(1).
\]
\end{theorem}

\begin{proof}
The proof consists of the following 3 steps:
\begin{enumerate}
    \item If $c_k < \frac{1}{4}$, then $c_{k+1} < \frac{1}{4}$.
    \item $c_k \to \frac{1}{4}$ as $k \to \infty$.
    \item If $\theta_k = \frac{k+2}{2} + \frac{\log k}{4} + e_k$, then $e_k$ is convergent.
\end{enumerate}
\paragraph{First step.}
If $c_k < \frac{1}{4}$, then $c_{k+1}<\frac{1}{4}$.

Let $\theta_k = \frac{k+2}{2} + c_k\log k $.
For our convenience, let $c_0=0$ with $c_0 \log 0 = 0$.
Plugging this in $\theta_{k+1}^2 - \theta_{k+1} - \theta_{k}^2=0$, we have
$$\left(\frac{k+2}{2} + c_{k+1} \log (k+1) \right)^2 = \left(\frac{k+2}{2} + c_{k}\log k \right)^2  + \frac{1}{4},$$
so
$$\left(c_{k+1} \log (k+1) - c_k \log k \right)\left(k+2+c_{k+1}\log (k+1) + c_k\log k \right) = \frac{1}{4}.$$
Assume $c_{k+1}\geq1/4$. Then 
\begin{align*}
  \frac{1}{4} &= \left(c_{k+1}\log(k+1) - c_k \log k \right)\left(k+2+c_{k+1}\log (k+1) + c_k \log k \right)
  \\
  &\geq \frac{1}{4}\log \left(1+\frac{1}{k}\right) (k+2)
  \\
  &>\frac{1}{4},
\end{align*}
which proves the first claim.

\paragraph{Second step.}
$c_k \to \frac{1}{4}$ as $k \to \infty$.

Put $d_k = \frac{1}{4}-c_k$, then $0 < d_k \le \frac{1}{4}$.
\begin{align*}
  \frac{1}{4} &=\left(\frac{1}{4}\log\left(1+\frac{1}{k}\right) -d_{k+1}\log (k+1) + d_k \log k \right) \left(k+2+\frac{1}{4}\log k(k+1)  -d_{k+1}\log (k+1) - d_k \log k \right)
  \\
  & \leq \left(\frac{1}{4}\log\left(1+\frac{1}{k}\right) -d_{k+1}\log (k+1) + d_k \log k \right) \left(k+2+\frac{1}{2}\log (k+1) \right)
\end{align*}
Therefore
\begin{align*}
  d_{k+1} \log(k+1) - d_k \log k \leq \frac{1}{4} \log \left(1+\frac{1}{k}\right) - \frac{1}{4}\frac{1}{k+2+\frac{1}{2}\log(k+1)}.
\end{align*}
By talyor expansion, 
\begin{align*}
  d_{k+1} \log(k+1) - d_k \log k \leq \frac{1}{4} \left(\frac{3+2\log k}{2k^2} + \mathcal{O}\left(\frac{1}{k^2}\right)\right).
\end{align*}
So, By summing all the above inequality from 1 to $k$,
$$d_{k+1} \log (k+1) \leq C$$
so $d_{k+1} < \frac{C}{\log (k+1)}$.
In conclusion, as $k \rightarrow \infty$, $d_k \rightarrow 0$.

\paragraph{Third step.}
If $\theta_k = \frac{k+2}{2} + \frac{\log k}{4} + e_k$, then, $e_k$ converges.

From the previous claim, we can say that for some sufficiently large $k$, $|e_k| < \frac{1}{6}\log k$.  
$$\left(\frac{k+2}{2} + \frac{1}{4}\log (k+1) + e_{k+1} \right)^2 = \left(\frac{k+2}{2} + \frac{1}{4} \log k + e_{k} \right)^2  + \frac{1}{4}$$
Then, 
\begin{align*}
  \frac{1}{4} &=\left(\frac{1}{4}\log\left(1+\frac{1}{k}\right) +e_{k+1} - e_k \right) \left(k+2+\frac{1}{4}\log k(k+1)   + e_{k+1} + e_k \right)
  \\
  &\leq \left(\frac{1}{4}\log\left(1+\frac{1}{k}\right) +e_{k+1} - e_k \right) \left(k+2+\frac{5}{6}\log (k+1) \right).
\end{align*}
So, $$e_{k+1} - e_k \geq \frac{1}{4\left(k+2+ \frac{5}{6}\log (k+1) \right) } - \frac{1}{4}\log \left( 1+ \frac{1}{k} \right) = -\frac{\frac{5}{6}\log k + \frac{3}{2}}{k^2} + \mathcal{O}\left(\frac{1}{k^2}\right). $$ 
Summing this for $k=1,\dots,k$, we get that $e_{k+1} > D$ for some constant $D$.
Moreover, 
\begin{align*}
  \frac{1}{4} &=\left(\frac{1}{4}\log\left(1+\frac{1}{k}\right) +e_{k+1} - e_k \right) \left(k+2+\frac{1}{4}\log k(k+1)   + e_{k+1} + e_k \right)
  \\
  &\geq \left(\frac{1}{4}\log\left(1+\frac{1}{k}\right) +e_{k+1} - e_k \right) \left(k+2\right)>\frac{1}{4} + (k+2)(e_{k+1} - e_k),
\end{align*}
which indicates that $e_{k+1} < e_k$. Since $\{e_k\}_{k=0}^\infty$ is a decreasing sequence with a lower bound, is converges.
\qed\end{proof}

\end{document}